\documentclass[11pt]{amsart}
\usepackage{amsfonts,latexsym}
\usepackage{amsmath}
\usepackage{amscd}
\usepackage{float,amsmath,amssymb,mathrsfs,bm,multirow,graphics}
\usepackage[dvips]{graphicx}
\usepackage[all]{xy}

\newcommand{\nequation}{\setcounter{equation}{0}}
\renewcommand{\theequation}{\mbox{\arabic{section}.\arabic{equation}}}
\newcommand{\R}{{\Bbb R}}

\newcommand{\C}{{\Bbb C}}
\newcommand{\Z}{{\Bbb Z}}

\newcommand{\proofbegin}{\noindent{\it Proof.\,\,}}
\newcommand{\proofend}{\hfill$\Box$\bigskip}

\newtheorem{theorem}{Theorem}[section]
\newtheorem{proposition}[theorem]{Proposition}
\newtheorem{lemma}[theorem]{Lemma}
\newtheorem{corollary}[theorem]{Corollary}

\newtheorem{remark}[theorem]{Remark}

\newtheorem{figuretext}{Figure}

\input epsf
\date{\today}
\title[Integrable evolution equations and peakons]
{Integrable evolution equations on spaces of tensor densities and their peakon solutions}

\author{Jonatan Lenells}
\address{J.L.: Department of Applied Mathematics and Theoretical Physics, 
University of Cambridge, Cambridge CB3 0WA, UK}
\email{j.lenells@damtp.cam.ac.uk} 

\author{Gerard Misio\l ek}
\address{G.M.: Department of Mathematics, University of Notre Dame, IN 46556, USA}
\email{gmisiole@nd.edu} 

\author{Feride T\i\u{g}lay}
\address{F.T.: Department of Mathematics, University of New Orleans, Lake Front, New 
Orleans, Louisiana 70148, USA and
Section de Math\'ematiques, {\'E}cole Polytechnique F{\'e}d{\'e}rale de Lausanne,
CH--1015 Lausanne, Switzerland}
\email{feride.tiglay@epfl.ch}

\begin{document}
\begin{abstract} 
We study a family of equations defined on the space of tensor densities
of weight $\lambda$ on the circle and introduce two integrable PDE.
One of the equations turns out to be closely related to the inviscid Burgers
equation while the other has not been identified in any form before.
We present their Lax pair formulations and describe their bihamiltonian
structures. We prove local wellposedness of the corresponding Cauchy
problem and include results on blow-up as well as global existence of
solutions. Moreover, we construct ``peakon'' and ``multi-peakon'' solutions for all
$\lambda \neq 0,1$, and ``shock-peakons'' for $\lambda = 3$. We argue that there is a natural
geometric framework for these equations that includes other well-known
integrable equations and which is based on V. Arnold's approach to
Euler equations on Lie groups.
\end{abstract}
\maketitle


\tableofcontents
\section{Introduction} 
\nequation

Integrability of an infinite-dimensional dynamical system typically 
manifests itself in several different ways such as 
the existence of a Lax pair formulation or a bihamiltonian structure, 
the presence of an infinite family of conserved quantities 
or at least the ability to write down explicitly some of its solutions. 
In this paper we introduce and study two 
nonlinear partial differential equations 
and show that they possess all the hallmarks of integrability 
mentioned above. 
The first of these equations we shall refer to as 
\emph{the $\mu$Burgers ($\mu$B) equation}\footnote{This equation is mentioned in Remark 3.9 of \cite{Holm-Stanley} and Remark 3.2 of \cite{lu} 
as the high-frequency limit of the Degasperis-Procesi equation, 
see \eqref{DP} below. Our terminology will be explained in Section \ref{sec:3HS}.}
\begin{equation} \tag{$\mu$B} 
\mu(u_t) - u_{txx} - 3u_xu_{xx} - uu_{xxx} =0,
\end{equation} 
and the second as \emph{the $\mu$DP equation} 
\begin{equation} \tag{$\mu$DP} 
\mu(u_t) -u_{txx} + 3\mu(u)u_x - 3u_x u_{xx} - u u_{xxx} = 0,
\quad\; 
\mathrm{where} 
\;\;\; 
\mu(u) = \int_0^1 u \, dx.
\end{equation} 
Here $u(t, x)$ is a spatially periodic real-valued function of a time variable $t$ and a space variable $x \in S^1 \simeq [0,1)$.
Both of these equations belong to a larger family that also includes 
the Camassa-Holm equation \cite{C-H} (see also \cite{F-F}) 
\begin{equation*} \tag{CH} 
u_t - u_{txx} +3uu_x - 2u_xu_{xx} - uu_{xxx} =0, 
\end{equation*}
the Hunter-Saxton equation \cite{H-S} 
\begin{equation*} \tag{HS} 
u_{txx} + 2u_x u_{xx} + uu_{xxx} = 0, 
\end{equation*} 
the Degasperis-Procesi equation \cite{dp} 
\begin{equation*} \tag{DP} \label{DP} 
u_t - u_{txx} + 4uu_x - 3u_x u_{xx} - uu_{xxx} =0, 
\end{equation*} 
as well as the $\mu$-equation which was derived recently in \cite{K-L-M}  
(we will refer to it as the $\mu$CH equation in this 
paper)\footnote{In \cite{K-L-M} the authors refer to $\mu$CH 
as the $\mu$HS equation.} 
\begin{equation*} \tag{$\mu$CH} 
u_{txx} -2\mu(u)u_x + 2u_x u_{xx} + uu_{xxx} = 0, 
\end{equation*}
all of which are known to be integrable. 

One of the distinguishing features of the CH and DP equations which makes 
them attractive among integrable equations 
is the existence of so-called ``peakon'' solutions.
In fact, CH and DP are the cases $\lambda =2$ and $\lambda =3$, respectively, 
of the following family of equations
\begin{equation} \label{CHDPfamily} 
u_t - u_{txx} + (\lambda + 1)uu_x = \lambda u_xu_{xx} + uu_{xxx}, \qquad \lambda \in \Z,
\end{equation}  
with each equation in the family admitting peakons (see \cite{dhh}) 
although only $\lambda =2$ and $\lambda =3$ are believed to be 
integrable (see \cite{dp}).
One of our results will show that each equation in the corresponding 
$\mu$-version of the family (\ref{CHDPfamily}) given by
\begin{equation}\label{eq:mulambda}  
\mu(u_t)  -u_{txx} + \lambda \mu(u)  u_x = \lambda u_x u_{xx} + uu_{xxx}, \qquad \lambda \in \Z,
\end{equation}
also admits peakon solutions.
The choices $\lambda = 2$ and $\lambda = 3$ yield the $\mu$CH and $\mu$DP equations, 
respectively, and as with (\ref{CHDPfamily}), we expect that these are 
the only integrable members of the family (\ref{eq:mulambda}). 
Moreover, we will show that the $\mu$DP equation admits shock-peakon solutions 
of a form similar to those known for DP, see \cite{lu}.


In Section \ref{Euler_Arnold}, we present a natural setting 
in which all the equations above can be formally described 
as evolution equations on the space of tensor densities 
(of different weight $\lambda$) over the Lie algebra of smooth vector fields 
on the circle. In Section \ref{Lax} we present Lax pairs for $\mu$DP and $\mu$B, 
establishing their integrability. 
In Section \ref{Hamilton} we describe the Hamiltonian structure of the equations 
in (\ref{eq:mulambda}). In particular, we consider the bihamiltonian structure of 
the $\mu$DP equation together with the associated infinite hierarchy 
of conservation laws.
In Section \ref{Cauchy} we study the periodic Cauchy problem of $\mu$DP; 
we prove local wellposedness in Sobolev spaces and show that 
while classical solutions of $\mu$DP break down for certain intial data, 
the equation admits global solutions for other data. 
In Section \ref{Waves} we construct multi-peakon solutions of (\ref{eq:mulambda}) 
as well as shock-peakons of $\mu$DP.
In Section \ref{sec:3HS} we discuss the $\mu$B equation and its properties 
and present a geometric construction to explain its close relation to 
the (inviscid) Burgers equation. 
Finally, in Section \ref{sec:EPDiff} we consider generalizations of the $\mu$CH 
and $\mu$DP equations to a multidimensional setting.

Left open is the question of physical significance, if any, of the two 
equations $\mu$DP and $\mu$B. It is possible that they may play a role in 
the mathematical theory of water waves (as CH, see \cite{C-H}) or find applications 
in the study of more complex equations (as e.g. HS, see \cite{Saxton-Tiglay}). 
We do not pursue these issues here.

Our approach draws heavily on \cite{K-M} and \cite{K-L-M} and 
the present work is in some sense a continuation of those two papers. 

\begin{figure}
\begin{center}
$$ 
\xymatrix{
\ar[dr] \text{CH}  & m_t = -um_x - 2u_x m &  \ar[dl]  \mu\text{CH} 
& \ar[dr] DP  & m_t = -um_x - 3u_x m & \ar[dl]  \mu\text{DP} & 
  	\\
& \text{HS} &&& 3\text{HS} & 
}
$$
\begin{figuretext}\label{trianglesfig} 
The $(\lambda=2)$ and $(\lambda=3)$ families of equations.
\end{figuretext}
\end{center}
\end{figure}
%


\section{A family of equations on the space of $\lambda$ densities} \nequation
\label{Euler_Arnold}

Perhaps the simplest way to introduce the equations that are 
the main object of our investigation is by analogy with the known cases. 
We shall therefore first briefly review Arnold's approach to 
the Euler equations on Lie groups 
and then describe a more general set-up intended to capture 
the $\mu$DP and the $\mu$B equations. 

%
%
%
%
%

In a pioneering paper Arnold \cite{A} presented a general framework 
within which it is possible to employ geometric and Lie theoretic techniques 
to study a variety of equations (ODE as well as PDE) 
of interest in mathematical physics 
(see also more recent expositions in \cite{K-M} or \cite{kw}). 
Arnold's principal examples were the equations of motion 
of a rigid body in $\mathbb{R}^3$ 
and the equations of ideal hydrodynamics. 

The formal set-up is the following. 
Consider a possibly infinite-dimensional Lie group $G$ 
(the ``configuration space'' of a physical system) 
with Lie algebra $\mathfrak{g}$. 
Choose an inner product 
$\langle \cdot, \cdot \rangle$ on $\mathfrak{g}$ 
(essentially, the ``kinetic energy'' of the system) 
and, using right- (or left-) translations, endow $G$ with 
the associated right- (or left-) invariant Riemannian metric. 
The motions of the system can now be studied either through 
the geodesic equation defined by the metric on $G$ 
(equivalently, the geodesic flow on its tangent or cotangent bundles) 
or else directly on the Lie algebra $\mathfrak{g}$ using 
Hamiltonian reduction.\footnote{We will be making use of both structures 
in this paper.} 
The equation that one obtains by this procedure on $\mathfrak{g}$ 
is called the \emph{Euler} (or \emph{Euler-Arnold}) \emph{equation}. 
Using the inner product it can be reformulated as an equation 
on the dual algebra $\mathfrak{g}^\ast$ as follows. 

Let $( \cdot, \cdot )$ be the natural pairing between 
$\mathfrak{g}$ and $\mathfrak{g}^\ast$ 
and let $A: \mathfrak{g} \to \mathfrak{g}^\ast$ denote 
the associated inertia operator\footnote{In infinite-dimensional examples 
$A$ is often some positive-definite self-adjoint pseudodifferential 
operator acting on the space of smooth tensor fields on a compact manifold.} 
determined by the formula 
$(Au,v) = \langle u, v \rangle$ 
for any $u, v \in \mathfrak g$. 
The Euler equation on $\mathfrak{g}^\ast$ reads 
\begin{equation*} \label{eq:ad}  \tag{E} 
m_t = - \mathrm{ad}^\ast_{A^{-1}m}m, 
\qquad 
m = A u \in \mathfrak{g}^\ast, 
\end{equation*} 
where $\mathrm{ad}^\ast: \mathfrak{g} \to \mathrm{End}(\mathfrak{g}^\ast)$ 
is the coadjoint representation of $\mathfrak{g}$ given by 
\begin{equation}  \label{ad*g}
(\mathrm{ad}^\ast_u m, v) = - (m, [u,v]) 
\end{equation}
for any $u, v \in \mathfrak{g}$ and $m \in \mathfrak{g}^\ast$. 
$\mathrm{ad}^\ast_u$ is the infinitesimal version of the coadjoint action 
$
\mathrm{Ad}^\ast: G \times \mathfrak{g}^\ast \to \mathfrak{g}^\ast
$
of the group $G$ on the dual algebra $\mathfrak{g}^\ast$. 

In our case $G$ will be the group $\mathrm{Diff}(S^1)$ 
of orientation-preserving diffeomorphisms of the circle 
whose Lie algebra $\mathfrak{g}$ 
is the space of smooth vector fields $T_e\mathrm{Diff}(S^1) = \mathrm{vect}(S^1)$. 
The dual $\mathrm{vect}^\ast(S^1)$ is the space of distributions on $S^1$ 
but we shall consider only its ``regular part'' which can be 
identified with the space of quadratic differentials 
$\mathcal{F}_2 = \left\{ m(x) dx^2: m \in C^\infty(S^1) \right\}$ 
with the pairing given by 
\begin{equation*} \label{pair2} 
\big( m  dx^2, v \partial_x \big) 
= 
\int_0^1 m(x) v(x) \, dx 
\end{equation*} 
(see e.g. Kirillov \cite{k}). 
The coadjoint representation of $\mathrm{vect}(S^1)$ on the regular part 
of its dual space is in this case precisely the action of $\mathrm{vect}(S^1)$ 
on the space of quadratic differentials. 
By a direct calculationm using \eqref{ad*g} we have 
\begin{equation}  \label{ad*} 
\mathrm{ad}^\ast_{u\partial_x}m dx^2 
= 
\left( u m_x + 2 u_x m \right) dx^2
\end{equation}
and the Euler equation \eqref{eq:ad} on $\mathfrak{g}^*$ 
%
%
takes the form 
\begin{equation} \label{CHtriple} 
m_t = -\mathrm{ad}^\ast_{A^{-1}m}m = -um_x - 2u_x m, 
\qquad 
m = Au. 
\end{equation} 

Equivalently, when rewritten on $\mathrm{vect}(S^1)$, equation \eqref{CHtriple} 
becomes 
\begin{equation} \label{EF2} 
Au_t + 2u_x Au + u A u_x =0 
\end{equation} 
so that with an appropriate choice of the inertia operator $A$ 
(which is equivalent to picking an inner product on $T_e\mathrm{Diff}$) 
\begin{equation}\label{Adef}
A = 
\begin{cases} 
1 - \partial_x^2 & \text{for CH}, 
  	\\  
\mu - \partial_x^2 & \text{for $\mu$CH}, 
	\\  
- \partial_x^2 & \text{for HS}, 
\end{cases}
\end{equation}
one recovers the CH, $\mu$CH and 
HS equations\footnote{Note that in the case of the HS equation 
the inertia operator is degenerate, 
see \cite{K-M}, Section 4, for details.}  
as Euler equations on $\mathrm{vect}(S^1)$; 
these results can be found in  \cite{Mis98}, \cite{K-M} and \cite{K-L-M}.

\smallskip 
We now want to extend this formalism to include the DP equation 
as well as the two equations of principal interest in this paper: 
$\mu$DP and $\mu$B. 
Admittedly, our construction 
does not have the same beautiful geometric interpretation as that of Arnold's 
and hence perhaps is not completely satisfactory. 

The main point is to replace the space of quadratic differentials 
with the space of all tensor densities on the circle 
of weight $\lambda$.\footnote{We refer to \cite{o}, \cite{ot} or \cite{gr} 
for basic facts about the space of tensor densities.} 
Recall that 
a tensor density of weight $\lambda\geq 0$ (respectively $\lambda<0$) 
on $S^{1}$ is a section of the bundle $\bigotimes^{\lambda} T^{*}S^{1}$ 
(respectively $\bigotimes^{-\lambda}TS^{1}$) and set 
\begin{equation*} 
\mathcal{F}_{\lambda} 
=
\left\{ m(x)dx^{\lambda}: m(x)\in C^{\infty}(S^{1}) \right\}.
\end{equation*} 
There is a well-defined action of the diffeomorphism group 
$\mathrm{Diff}(S^1)$ on each density module $\mathcal{F}_\lambda$ 
given by 
\begin{equation}  \label{3-action} 
\mathcal{F}_\lambda \ni m dx^\lambda 
\to 
m\circ\xi \, (\partial_x\xi)^\lambda dx^\lambda \in \mathcal{F}_\lambda, 
\qquad 
\xi \in \mathrm{Diff}(S^1), 
\end{equation}
which naturally generalizes the coadjoint action 
$
\mathrm{Ad}^\ast: 
\mathrm{Diff}(S^1) \to \mathrm{Aut}\big(\mathcal{F}_2 \big)
$ 
on the space of quadratic differentials. 
The infinitesimal generator of the action in \eqref{3-action} 
is easily calculated,
\begin{equation} \label{Lie*} 
L_{u\partial_x}^{\lambda}(m dx^{\lambda}) 
= 
\left( um_{x}+\lambda u_{x}m \right) dx^{\lambda}, 
\end{equation} 
and can be thought of as the Lie derivative of tensor densities. 
It represents the action of $\mathrm{vect}(S^{1})$ on $\mathcal{F}_{\lambda}$ 
which for $\lambda=2$ coincides with the (algebra) coadjoint action 
on $\mathcal{F}_2$ 
(that is $L^2_{u\partial_x}=\mathrm{ad}^\ast_{u\partial_x}$). 

If we think of \eqref{Lie*} as defining a vector field 
on the space $\mathcal{F}_\lambda$ then we can consider 
the equation for its flow
\begin{equation} \label{EF3} 
m_t  = - um_{x} - \lambda u_{x}m 
\end{equation} 
in analogy with \eqref{CHtriple}. 
The substitution $m=Au$ transforms (\ref{EF3}) into an equation 
on the space of quadratic differentials 
\begin{equation} 
Au_t + \lambda u_x A u + u A u_x =0 
\label{eq:lambda}
\end{equation}
which is the $\lambda$-version of \eqref{EF2}. 

Setting $\lambda=3$ and choosing suitable inertia operators $A$ as above 
we obtain the DP, $\mu$DP and $\mu$B equations. 
More precisely, the specific choices of $A$ are in parallel 
with those made for the $\lambda=2$ family in \eqref{Adef}, that is 
\begin{equation}\label{Adeff}
A = 
\begin{cases} 
1 - \partial_x^2 & \text{for DP},
  	\\  
\mu - \partial_x^2 & \text{for $\mu$DP},
	\\  
- \partial_x^2 & \text{for $\mu$B}.
\end{cases}
\end{equation}

More generally, letting $A = \mu - \partial_x^2$ in (\ref{eq:lambda}) for any $\lambda \in \Z$, we find the equations in (\ref{eq:mulambda}).


\section{Lax pairs}\nequation
\label{Lax}

An elegant manifestation of complete integrability 
of an infinite dimensional dynamical system is the existence of 
a Lax pair formalism. 
It is often used as a tool for constructing infinite families of 
conserved quantities. 
This formalism has been quite extensively developed in recent years 
for the CH and DP equations, see e.g. \cite{bss, cm, dhh, ls}. 
Our next result describes the Lax pair formulations for the $\mu$DP and 
$\mu$Burgers equations. 

\begin{theorem}  \label{thm:Lax} 
The $\mu$DP and the $\mu$B equations admit the Lax pair formulations 
\begin{equation}\label{lax}
\begin{cases} 
\psi_{xxx}= - \lambda m \psi, 
\\
\psi_t 
= 
-\frac{1}{\lambda} \psi_{xx} - u\psi_x + u_x \psi, 
\end{cases}
\end{equation}
where $\lambda \in \C$ is a spectral parameter, $\psi(t, x)$ is a scalar eigenfunction 
and $m = Au$ with $A=\mu -\partial_x^2$ or $A=-\partial_x^2$ as defined in \eqref{Adeff}.
\end{theorem} 
\proofbegin
This is a straightforward computation which shows that the condition of compatibility 
$(\psi_t)_{xxx} = (\psi_{xxx})_t$ of the linear system (\ref{lax}) is equivalent to 
the equation $\mu$DP or $\mu$B when $m$ is given by 
$m=\mu(u) - u_{xx}$ or $m=-u_{xx}$, respectively.
\proofend


\section{Hamiltonian structures and conserved quantities}\nequation
\label{Hamilton}
In this section we describe the bihamiltonian structure of $\mu$DP 
(see Section \ref{sec:3HS} for the bihamiltonian structure of 
the $\mu$Burgers equation) 
as well as the Hamiltonian structure for the family of 
$\mu$-equations (\ref{eq:mulambda}) for any $\lambda \in \mathbb{Z}$. 
In the last subsection we consider the orbits of $\mathrm{Diff}(S^1)$ 
in the space of tensor densities of weight $\lambda$ and give a geometric interpretation of 
one of the conserved quantities.


\subsection{Bihamiltonian structure of $\mu$DP}   
Recall that DP admits the bihamiltonian formulation \cite{dhh}
$$
m_t = J_0 \frac{\delta H_0}{\delta m} = J_2 \frac{\delta H_2}{\delta m},
$$
where 
$$
H_0 
= 
-\frac{9}{2} \int m \, dx 
\qquad 
\mathrm{and}
\qquad 
H_2 
= 
-\frac{1}{6} \int u^3 dx, 
$$
and the Hamiltonian operators are
$$
J_0 
= 
m^{2/3}\partial_x m^{1/3}
\big( \partial_x - \partial_x^3 \big)^{-1} m^{1/3}\partial_x m^{2/3} 
\quad
\mathrm{and}
\quad 
J_2 = \partial_x \big( 4 - \partial_x^2 \big)\big(1 - \partial_x^2 \big) 
$$
with $m = (1 - \partial_x^2)u$. 
\begin{figure}
\begin{center}
$$ \xymatrix{
 \mu\text{DP}_{3}  &   
\\
  &  \ar@/^/[dl]^{\qquad \qquad \  J_0}  \frac{\delta  H_2}{\delta m}	
  		\\
m_t = -m_x u - 3 m u_x & & 
 		\\  
    &   \ar@/_/[uuul]^{J_2 \qquad}  \ar@/^/[dl]^{\qquad \qquad \ J_0}  \frac{\delta}{\delta m}\left( \frac{1}{2} \int u^2 dx\right) = \frac{\delta H_1}{\delta m}\\
m_t = m_x &  
 		\\
  &    \ar@/_/[uuul]^{J_2 \qquad} \ar@/^/[dl]^{\qquad \qquad \ J_0}   \frac{\delta}{\delta m}\left( -\frac{9}{2} \int m dx\right) = \frac{\delta H_0}{\delta m}
   		\\
m_t = 0   & 
		 \\
 &    \ar@/_/[uuul]^{J_2 \qquad} \ar@/^/[dl]^{\qquad \qquad \ J_0}    \frac{\delta}{\delta m}\, \int m^{1/3} dx = \frac{\delta H_{-1}}{\delta m}	
		 \\
\mu\text{DP}_{-1}   &   
 		 \\
 &    \ar@/_/[uuul]^{J_2 \qquad}  \frac{\delta}{\delta m}\left(-\frac{1}{2}\int \frac{m_x^2}{m^{7/3}} dx\right) = \frac{\delta H_{-2}}{\delta m}
 		 }$$
     \begin{figuretext}\label{2hierarchyfig}
       Recursion scheme for the $\mu$DP equation. 
     \end{figuretext}
     \end{center}
\end{figure}

Similarly, the $\mu$DP equation admits the bihamiltonian formulation
$$
m_t 
= 
J_0 \frac{\delta H_0}{\delta m} 
= 
J_2 \frac{\delta H_2}{\delta m},
$$
where now the Hamiltonian functionals $H_0$ and $H_2$ are 
\begin{equation}\label{H0H2muDPdef}
H_0 
= 
-\frac{9}{2} \int m \, dx 
\qquad 
\mathrm{and} 
\qquad 
H_2 
=  
- \int \left(
\frac{3}{2} \mu(u) \big(A^{-1}\partial_x u\big)^2 
+ 
\frac{1}{6}u^3 
\right)dx,
\end{equation}
the operators $J_0$ and $J_2$ are given by
$$
J_0 
= -m^{2/3}\partial_x m^{1/3}\partial_x^{-3} m^{1/3}\partial_x m^{2/3} 
\qquad
\mathrm{and} 
\qquad 
J_2 = - \partial_x^3 A = \partial_x^5 
$$
and $m=Au$ with $A = \mu - \partial_x^2$. 
The fact that $J_0$ and $J_2$ form a compatible bihamiltonian pair 
is a consequence of Theorem 2 in \cite{hw}.
It can be verified directly that $H_0$ and $H_2$, as well as 
$H_1 = \frac{1}{2} \int u^2 dx$
are conserved in time whenever $u$ is a solution of $\mu$DP (see Appendix \ref{H2app} for details of this calculation in the case of $H_2$).
Using the standard techniques we can now construct an infinite sequence of 
conservation laws 
$$
\dots H_{-1}, \, H_0, \, H_1, \, H_2, \dots 
$$ 
see Figure \ref{2hierarchyfig}.
As in the case of the CH and DP equations the above conservation laws $H_n$ 
are nonlocal and not easy to write down explicitly for $n \geq 3$, 
while they are readily computable recursively in terms of $m$ and its derivatives 
for negative $n$. 
In fact, the first negative flow in the $\mu$DP hierarchy is
\begin{align*}
\mu\text{DP}_{-1}: m_t 
=& 
-\frac{2}{729 m^{17/3}} \biggl(6160 m_x^5-13200 m m_{xx} m_x^3+3600 m^2 m_x^2 m_{xxx}
	\\
&-675 m^2 \left(m m_{xxxx} - 8  m_{xx}^2\right) m_x
+27 m^3 \left(3 m m_{xxxxx} - 50 m_{xx} m_{xxx} \right)\biggr).
\end{align*}
%


\subsection{Hamiltonian formulation of $\mu$-equations}
Although we expect that $\mu$CH and $\mu$DP are the only equations 
among the family of $\mu$-equations in (\ref{eq:mulambda})
which admit a bihamiltonian structure 
we can nevertheless provide one Hamiltonian structure for any $\lambda \neq 1$. 
Indeed, if we set
\begin{equation} \label{J0lambdadef}
J_0 
= 
-\frac{1}{\lambda^2}(m_x + \lambda m\partial_x) 
\partial_x^{-3}\big( 
(\lambda-1)m_x + \lambda m\partial_x 
\big) 
\end{equation}
and
\begin{equation} 
H_0 = -\frac{\lambda^2}{\lambda-1} \int m \, dx 
\end{equation}
then equation (\ref{eq:mulambda}) 
is equivalent to
$$
m_t = J_0 \frac{\delta H_0}{\delta m}.
$$
In order to see that $J_0$ is a Hamiltonian operator we first rewrite it 
in the form 
$$
J_0 
= 
-m^{(\lambda-1)/\lambda}\partial_x m^{1/\lambda}\partial_x^{-3}m^{1/\lambda}
\partial_x m^{(\lambda-1)/\lambda}
$$
and refer to Theorem 1 of \cite{hw}.


\subsection{Orbits in $\mathcal{F}_\lambda$}
From Section \ref{Euler_Arnold} we know that for any value of $\lambda$ 
a solution $m$ of 
the Cauchy problem\footnote{See Section \ref{subsec:lwp} for results 
on wellposedness of the equations \eqref{eq:mulambda} for any $\lambda$.} 
for equation (\ref{eq:mulambda}) with initial data $m_0$ belongs to the orbit of 
$m_0$ in $\mathcal{F}_\lambda$ under the $\mathrm{Diff}(S^1)$-action 
defined in (\ref{3-action}). 
It is therefore tempting to exploit this setting to study the family 
in \eqref{eq:mulambda} as in the case of the coadjoint action when $\lambda =2$. 
We present here two results in this direction leaving a detailed investigation 
of the geometry of the corresponding orbits for general $\lambda$ for a future work.

The first proposition shows the conservation of 
$H_{-1} = \int_0^1 |m|^{1/\lambda} dx$ 
under the flow of (\ref{eq:mulambda}). 
We let 
$L_\xi(mdx^\lambda) = m \circ \xi (\partial_x \xi)^\lambda dx^\lambda$ 
denote the action of an element 
$\xi \in \mathrm{Diff}(S^1)$ on $m dx^\lambda$ in $\mathcal{F}_\lambda$.
\begin{proposition}\label{orbitprop1}
The map
$$ 
m \, dx^\lambda \mapsto H_{-1}[m] 
= 
\int_0^1 |m|^{1/\lambda} dx:\mathcal{F}_\lambda \to \R_+ 
$$ 
is invariant under the action (\ref{3-action}) of $\mathrm{Diff}(S^1)$ 
on $\mathcal{F}_\lambda$.
\end{proposition}
\proofbegin 
This is a straightforward change of variables 
$$
H_{-1} \big[L_\xi (m \, dx^\lambda) \big] 
= 
\int_0^1 \left| m \circ \xi \, (\partial_x \xi)^\lambda \right|^{1/\lambda} dx 
= 
\int_0^1 |m|^{1/\lambda} \circ \xi \, d\xi  
= 
H_{-1}[m],
$$ 
since $\xi$ is a smooth orientation-preserving circle diffeomorphism 
and $m \, dx^\lambda \in \mathcal{F}_\lambda$.
\proofend

In many respects the classification of orbits in $\mathcal{F}_\lambda$ 
resembles that of coadjoint orbits of $\mathrm{Diff}(S^1)$, 
cf. \cite{gr}. 
In order to state the next result 
we denote by $\mathcal{M}_\lambda$ the set of those elements 
$m dx^\lambda$ in $\mathcal{F}_\lambda$ such that $m(x) \neq 0$ 
for all $x \in S^1$. 
\begin{proposition}\label{orbitprop2}
The orbit space $\mathcal{M}_\lambda / \mathrm{Diff}(S^1)$ is in bijection 
with the set $\R\backslash \{0\}$. 
More precisely, the map
\begin{equation}\label{orbitmap}
m dx^\lambda \mapsto \mathrm{sgn}(m) (H_{-1})^\lambda: 
\mathcal{M}_\lambda \to \R\backslash \{0\}, 
\qquad 
H_{-1} = \int_0^1 |m|^{1/\lambda} dx,
\end{equation}
is constant on the $\mathrm{Diff}(S^1)$-orbits in 
$\mathcal{M}_\lambda \subset \mathcal{F}_\lambda$ 
and induces a one-to-one correspondence between orbits in 
$\mathcal{M}_\lambda$ and $\R\backslash \{0\}$.
\end{proposition}
\proofbegin 
That the map is constant on the orbits follows from Proposition \ref{orbitprop1}. 
To see that it is surjective we note that for any nonzero real number 
$a$ the element $a dx^\lambda$ is mapped to $a$. 
In order to prove injectivity we observe that the orbit through an arbitrary element 
$m dx^\lambda$ of $\mathcal{M}_\lambda$ contains an element of the form 
$a dx^\lambda$ where $a$ is a constant. 
Indeed, given $m dx^\lambda$ we define $\xi \in \mathrm{Diff}(S^1)$ 
by
$$
\xi(x) = \frac{1}{H_{-1}} \int_0^x |m|^{1/\lambda} dx.
$$
Then
$$
L_\xi\left(\mathrm{sgn}(m) H_{-1}^\lambda dx^\lambda\right) 
= 
\mathrm{sgn}(m) H_{-1}^\lambda (\partial_x \xi)^\lambda dx^\lambda 
= 
m dx^\lambda.
$$
Thus, the action of $\xi^{-1}$ maps $m$ to 
$\mathrm{sgn}(m) H_{-1}^\lambda dx^\lambda$ 
where $\mathrm{sgn}(m) H_{-1}^\lambda \in \R\backslash\{0\}$ is constant.
\proofend

The stabilizer in $\mathrm{Diff}(S^1)$ of a constant element 
$a dx^\lambda$ in $\mathcal{M}_\lambda$ $(a \in \R)$ consists exactly of 
the rigid rotations and is isomorphic to $S^1$.
It follows from Proposition \ref{orbitprop2} that the stabilizer of any point 
$m dx^\lambda\in \mathcal{M}_\lambda$ is conjugate to $S^1$ within 
$\mathrm{Diff}(S^1)$. 
In particular, each orbit in $\mathcal{M}_\lambda$ is of the form 
$\mathrm{Diff}(S^1)/S^1$.


\section{The Cauchy problem for the $\mu$DP equation}  \nequation
\label{Cauchy} 

In this section we turn to the periodic Cauchy problem for $\mu$DP. 
Rather than aiming at the strongest possible theorems we present here 
only basic results that display the interesting behavior of solutions. 

We start with local existence, uniqueness and persistence theorems 
for the $\mu$DP equation in Sobolev spaces. 
Next, we describe the breakdown of smooth solutions and 
prove a global existence result for the class of initial data 
with non-negative momentum density. 
At the end of this section we briefly discuss local wellposedness 
of the whole family (\ref{eq:mulambda}) of $\mu$-equations. 

In order to prove these results we will need Sobolev completions 
of the group of 
circle diffeomorphisms\footnote{Although no longer a Lie group, 
$\mathrm{Diff}^s(S^1)$ retains the structure of a topological group 
for a sufficiently high Sobolev index; see below.} 
 $\mathrm{Diff}(S^1)$
and its Lie algebra of smooth vector fields $\mathrm{vect}(S^1)$. 
We denote by $H^s=H^s(S^1)$ the Sobolev space of periodic functions 
\begin{equation*} 
H^s(S^1) 
= 
\Big\{ v= \sum_n \hat{v}(n) e^{2\pi i nx}: 
\|v\|_{H^s}^2 = \sum_n \big| \widehat{\Lambda^s v} (n) \big|^2 < \infty 
\Big\}
\end{equation*}
where $s \in \mathbb{R}$ and  the pseudodifferential operator 
$\Lambda^s = (1 - \partial_x^2)^{s/2}$ 
is defined by 
$$
\widehat{\Lambda^s v}(n) = ( 1 + 4\pi^2n^2 )^{s/2} \hat{v}(n). 
$$
In what follows we will also make use of another elliptic operator 
\begin{equation} \label{lambda-mu}
\Lambda_\mu^2 : H^s(S^1) \to H^{s-2}(S^1), 
\qquad 
\Lambda_\mu^2 v = \mu(v) - v_{xx} 
\end{equation}
whose inverse can be easily checked to be 
\begin{align} \label{inverse} 
\Lambda_\mu^{-2}v(x) 
= 
\left( \frac{x^2}{2} -\frac{x}{2} +\frac{13}{12} \right)\int_0^1 v(x)\, dx 
+ 
\Big( x &-\frac{1}{2} \Big) \int_0^1 \int_0^x v(y)\, dy dx     \\ \nonumber 
- 
\int_0^x \int_0^y v(z) \, dz dy 
&+ 
\int_0^1 \int_0^x \int_0^y v(z) \, dz dy dx.  
\end{align} 
%


\subsection{Local wellposedess in Sobolev spaces} 

We begin with the $\mu$DP equation. Our proof follows with minor changes 
the approach of \cite{Mis02} to the CH equation based on the original methods 
in \cite{em} 
developed for the Euler equations of hydrodynamics. 
We include it for the sake of completeness. 
The two above papers and \cite{K-L-M} will be our main references 
for most of the basic facts we use in this section. 

We can write the Cauchy problem for $\mu$DP in the form
\begin{align} \label{muDP} 
u_t + u u_x + 3\mu(u)\, \partial_x \Lambda_\mu^{-2} u = 0    
\end{align} 
\begin{equation}  \label{ic} 
u(0)=u_0. 
\end{equation} 
Note that applying $\Lambda_\mu^2$ to both sides of \eqref{muDP} we obtain 
the equation in its original form given in the Introduction. 

\begin{theorem}[Local wellposedness and persistence] 
Assume $s > 3/2$. Then for any $u_0\in H^s(\mathbb{T})$ there exists a $T > 0$ 
and a unique solution 
$$
u\in C\big( (-T , T ), H^s \big) \cap C ^1\big( ( -T , T ), H^{s-1}\big) 
$$
of the Cauchy problem \eqref{muDP}-\eqref{ic} which depends continuously 
on the initial data $u_0$. 
Furthermore, the solution persists as long as $\|u (t , \cdot )\|_{C^ 1}$ 
stays bounded. 
\label{th:lwp}
\end{theorem}

Our strategy will be to reformulate \eqref{muDP}-\eqref{ic} 
as an initial value problem on the space of circle diffeomorphisms $\mathrm{Diff}^s(S^1)$ 
of Sobolev class $H^s$. 
It is well known that whenever $s>3/2$ this space is a smooth Hilbert manifold 
and a topological group. 
We will then show that the reformulated problem can be solved 
on $\mathrm{Diff}^s(S^1)$ by standard ODE techniques. 

Let $u=u(t,x)$ be a solution of $\mu$DP with initial data $u_0$. 
Then its associated flow $t \to \xi(t, x)$, i.e. the solution of 
the initial value 
problem\footnote{Here ``dot'' indicates differentiation in $t$ variable.} 
\begin{equation} \label{mu-flow} 
\dot{\xi}(t,x) = u(t, \xi(t,x)), 
\quad 
\xi(0,x) = x
\end{equation} 
is (at least for a short time) a $C^1$ smooth curve in the space of diffeomorphisms 
starting from the identity $e \in \mathrm{Diff}^s(S^1)$. 
Differentiating both sides of this equation in $t$ and using \eqref{muDP} 
we obtain the following initial value problem 
\begin{align}
& \ddot{\xi} 
= 
-3\mu \big(\dot{\xi} \circ \xi^{-1} \big) 
\Big( \partial_x \Lambda_{\mu}^{-2} (\dot{\xi} \circ \xi^{-1}) \Big)\circ\xi 
=: -F(\xi,\dot\xi)  
\label{eq:ode} \\
& \xi(0,x)=x, \quad  \dot{\xi}(0,x)=u_{0}(x)  
\label{iic}
\end{align}

In this section we make repeated use of a technical result which we state here 
for convenience (for a proof refer to \cite{Mis02}, Appendix 1). 
\begin{lemma} For $s>3/2$ the composition map $\xi\rightarrow \omega\circ\xi$ 
with an $H^s$ function $\omega$ and the inversion map $\xi\rightarrow \xi^{-1}$ 
are continuous from $\mathrm{Diff}^s(S^1)$ to $H^s(S^1)$ and from $\mathrm{Diff}^s(S^1)$  
to itself respectively. Moreover,  
\begin{equation}
\| \omega\circ\xi \|_{H^s}\leq C(1+\|\xi\|_{H^s}^s)\|\omega\|_{H^s}
\end{equation}
with $C$ depending only on $\inf |\partial_x \xi |$ and $\sup |\partial_x \xi| $.
\label{l:comp}
\end{lemma}

It will also be convenient to introduce the notation 
$A_{\xi}=R_{\xi}\circ A \circ R_{\xi^{-1}}$
for the conjugation of an operator $A$ on $H^s(S^1)$ by 
a diffeomorphism $\xi \in \mathrm{Diff}^s(S^1)$. 
The right hand side of \eqref{eq:ode} then becomes 
\begin{equation}
F(\xi, \dot{\xi})=\Lambda^{-2}_{\mu,\xi} \partial_{x,\xi} h(\xi,\dot{\xi})
\label{eq:F}
\end{equation}
where $h(\xi,\omega)=3\omega \int_0^1 \omega \circ \xi^{-1} dx$. 

\medskip
\noindent{\it Proof of Theorem \ref{th:lwp}.} 
From Lemma \ref{l:comp} we readily see that $F$ maps into $H^s(S^1)$. 
We aim to prove that $F$ is Fr\'{e}chet differentiable in 
a neighborhood of the point $(e,0)$ in $\mathrm{Diff}^{s}\times H^{s}(S^1)$.  
To this end we compute the directional derivatives 
$\partial_{\xi} F_{(\xi,\omega)}$ and $ \partial_{\omega} F_{(\xi,\omega)}$
and show that they are bounded linear operators on $H^s$ which depend continuously 
on $\xi$ and $\omega$. 
We use the formulas
\begin{equation}
\partial_{\xi} \Lambda^{-2}_{\mu,\xi}(v) 
= 
-\Lambda^{-2}_{\mu,\xi}
\left[ v\circ\xi^{-1}\partial_x , \Lambda_\mu^2 \right]_{\xi}
\Lambda^{-2}_{\mu,\xi} 
\end{equation}
\begin{equation}
\partial_{\xi} \partial_{x,\xi}(v) 
= 
\left[ v\circ\xi^{-1}\partial_x, \partial_x \right]_{\xi}
\end{equation}
and
\begin{equation}
\partial_{\xi} h_{(\xi,\omega)}(v) 
= 
3\, \omega \int\omega\circ\xi^{-1} \partial_x (v\circ\xi^{-1}) \, dx
\end{equation}
along with \eqref{eq:F} to get 
\begin{align*} 
\partial_{\xi} F_{(\xi,\omega)}(v) 
=
3 \, \Big\{ 
&-v\circ\xi^{-1}\Lambda^{-2}_\mu\partial_x^2(\omega \circ \xi^{-1}) 
\int \omega\circ\xi^{-1} dx        \\
&+ \Lambda^{-2}_\mu \partial_x \big( 
(v\circ\xi^{-1})\partial_x (\omega\circ\xi^{-1}) \big) 
\int (\omega\circ\xi^{-1})\, dx                   \\ 
&- \Lambda^{-2}_\mu 
\partial_x (\omega\circ\xi^{-1})
\int(\omega\circ\xi^{-1})\, \partial_x (v\circ\xi^{-1})\, dx 
\Big\}\circ\xi 
\end{align*}
which with the help of the identity
\begin{equation}  \label{A}
\Lambda^{-2}_\mu \partial_x^2 =  -1+\mu 
\end{equation}
can be simplified to 
\begin{align*}
\partial_{\xi} F_{(\xi,\omega)}(v)
=
3 \left\{ 
v \,\omega \int \omega\circ\xi^{-1} dx                             \right. 
&- 
v \left( \int \omega\circ\xi^{-1} dx \right)^2      \\
&+  
\Lambda^{-2}_{\mu,\xi} \partial_{x,\xi} \big( v \, \partial_{x,\xi} \omega \big) 
\int \omega\circ\xi^{-1} dx           \\ 
&-                                                                  \left. 
\Lambda^{-2}_{\mu,\xi} \partial_{x,\xi} \omega 
\int \omega\circ\xi^{-1} \partial_x (v\circ\xi^{-1}) \, dx 
\right\}. 
\end{align*}

On the other hand we have
\begin{equation*} 
\partial_{\omega}h_{(\xi,\omega)}(v) 
= 
3 v \int \omega\circ\xi^{-1} dx + 3 \omega \int v\circ\xi^{-1} dx 
\end{equation*}
and similarly 
\begin{equation*}
\partial_{\omega} F_{(\xi,\omega)}(v) 
= 
-3 \Lambda^{-2}_{\mu,\xi} \partial_{x,\xi} v \int \omega\circ\xi^{-1} dx 
-
3 \Lambda^{-1}_{\mu,\xi} \partial_{x,\xi}\omega \int v\circ\xi^{-1} dx.
\end{equation*}

In order to show that $v \rightarrow \partial_{\xi} F_{(\xi ,\omega)} (v)$ 
is a bounded operator on $H^s$ it is sufficient to estimate the sum 
\begin{equation*}
\| \partial_{\xi}F_{(\xi,\omega)}(v) \|_{L^2} 
+ 
\| \partial_x\big(
\partial_{\xi}F_{(\xi,\omega)}(v)\big)\circ\xi^{-1} \|_{H^{s-1}}.
\end{equation*}
Using Cauchy-Schwarz, Lemma \ref{l:comp} and the formulas above 
both of these terms can be bounded by a multiple of
$
\|v\|_{H^s} \|\omega\|_{H^s}^2 
$
which gives the estimate 
$$
\| \partial_{\xi} F_{(\xi ,\omega)} (v)\|_{H^s} 
\leq 
C \|v\|_{H^s} \|\omega \|_{H^s}^2 
$$
where $C$ depends only on $\inf |\partial_x \xi |$ and $\sup |\partial_x \xi|$.
A similar argument yields the bound for the other directional derivative 
$$
\| \partial_{\omega}F_{(\xi,\omega)}(v)\|_{H^s} 
\leq 
C \|v\|_{H^s} \|\omega \|_{H^s}.
$$

It follows that $F$ is Gateaux differentiable near $(e,0)$ and therefore 
we only need to establish continuity of the directional derivatives. 
Continuity in the $\omega$-variable follows from the fact that the dependence 
of both partials on this variable is polynomial. 
It thus only remains to show that the norm of the difference 
\begin{equation*}
\big\| \partial_{\xi}F_{(\xi,\omega)}(v) 
-
\partial_{\xi}F_{(e,\omega)}(v) \big\|_{H^s}
\end{equation*}
is arbitrarily small whenever $\xi$ is close to the identity 
$e \in \mathrm{Diff}^s(S^1)$, uniformly in $v$ and $\omega \in H^s(S^1)$. 
Once again, it suffices to estimate the terms 
\begin{equation}
\big\| \partial_{\xi}F_{(\xi,\omega)}(v) 
- 
\partial_{\xi}F_{(e,\omega)}(v) \big\|_{L^2} 
+  
\big\| \partial_x\big(\partial_{\xi}F_{(\xi,\omega)}(v) 
- 
\partial_{\xi}F_{(e,\omega)}(v) \big) \big\|_{H^{s-1}}
\label{eq:dFcont}
\end{equation}
We proceed as above. Using formula \eqref{inverse} 
and Lemma \ref{l:comp} together with 
the fact that for any $r>0$ we have a continuous embedding 
$H^r ( S^1 ) \hookrightarrow C^{r-1/2}(S^1)$  
with a pointwise estimate 
$$
\|g(x ) - g( y )\|  \lesssim \| g\|_{H^r}\big| x-y \big|^{r-1/2}, 
$$
we can bound both terms in \eqref{eq:dFcont} by
$$
C \| \omega \|_{H^s}^2 \|v\|_{H^s} \left(
\big\| \xi^{-1} - e \big\|_{H^s} 
+  
\big\| \xi^{-1} - e \big\|_{\infty}^{s-3/2}
\right).
$$
Continuity of $\xi\rightarrow \partial_{\omega}F_{(\xi,\omega)}$ follows from 
a similar estimate.

In this way we obtain that $F$ is continuously differentiable near $(e, 0)$ 
and applying the fundamental ODE theorem for Banach spaces 
(see e.g. Lang \cite{l}) 
establish local existence, uniqueness and smooth dependence on the data $u_0$ 
of solutions $\xi(t)$ and $\dot{\xi}(t)$ of \eqref{eq:ode}-\eqref{iic}. 
Local wellposedness of the original Cauchy problem \eqref{muDP}-\eqref{ic} 
follows now at once from $u=\dot\xi\circ\xi^{-1}$ and the fact that 
$\mathrm{Diff}^s$ is a topological group whenever $s>3/2$. 


In order to complete the proof of Theorem \ref{th:lwp} we need to derive 
a $C^1$ bound on the solution $u$. 
The standard trick is to use Friedrichs' mollifiers $J_{\epsilon}$ 
with ($0<\epsilon <1$) to derive the inequality 
\begin{align*}
\frac{d}{dt} \left\| J_{\epsilon}u \right\|_{H^s}^2 
&= 
- \int_0^1\Lambda^s J_{\epsilon}u\,\partial_x\Lambda^s J_{\epsilon}u^2 dx    
- 
\frac{3}{2} \mu(u) \int_0^1 
\Lambda^2 J_\epsilon u\,\partial_x\Lambda^s\Lambda_\mu^{-2}J_\epsilon u 
\, dx                             \\
&\lesssim 
\|u\|_{C^1} \|u\|_{H^s}^2 
+ 
|\mu(u)| 
\left\|\partial_x\Lambda^s\Lambda_\mu^{-2}J_\epsilon u\right\|_{L^2} 
\left\|\Lambda^s J_\epsilon u \right\|_{L^2} 
\lesssim 
\|u\|_{C^1} \|u\|_{H^s}^2 
\end{align*}
where the first estimate is standard (see e.g. \cite{t}, Chapter 5) 
and the second follows easily with the help of 
\eqref{lambda-mu} and \eqref{A}. 
Passing to the limit with $\epsilon \to 0$ and 
using Gronwall's inequality yields the persistence result. 
\proofend


\subsection{Blow-up of smooth solutions}

A family of examples displaying a simple finite-time breakdown mechanism 
of a $\mu$DP solution is described by the following result.  
\begin{theorem} \label{muDPblowup} 
Given any smooth periodic function $u_0$ with zero mean there exists 
$T_c >0$ such that the corresponding solution of the $\mu$DP equation 
stays bounded for $t<T_c$ and satisfies 
$\| u_x(t)\|_\infty \nearrow \infty$ as $t \nearrow T_c$.
\end{theorem}
\proofbegin  Let $u$ be the solution of $\mu$DP with initial data $u_0$. 
Differentiating the equation with respect to the space variable we obtain
$$
u_{tx} + u u_{xx} + u_x^2 + 3 \mu(u) \partial_x^2 \Lambda_\mu^{-2} u = 0.
$$
Since 
$\partial_x^2 \Lambda_\mu^{-2} = \mu - 1$,
using conservation of the mean, we find that
\begin{equation} \label{1}
u_{tx} + u u_{xx} + u_x^2 = 3 \mu(u) \big( u - \mu(u) \big) = 0 
\end{equation} 
by the assumption on the data $u_0$. 
If we let $\xi(t)$ denote the flow of $u$ as in \eqref{mu-flow} 
then  
$$
\partial_x\dot{\xi} = u_x\circ\xi \, \partial_x \xi 
$$
and setting 
$w=\partial_x  \dot{\xi}/\partial_x \xi$ 
we find that 
$$
w_t 
= 
\frac{\partial_x\ddot\xi\,\partial_x\xi-(\partial_x\dot\xi)^2}{(\partial_x\xi)^2}
= 
\big( u_{xt} + u_{xx} u \big)\circ\xi. 
$$
With the help of these formulas we rewrite \eqref{1} in the form 
$$
w_t+w^2=0 
$$
and, since $w(0,x)=u_{0 x}(x)$, solve for $w$ to get 
$$
w(t,x)=\frac{1}{t+\big(1/u_{0x}(x)\big)}. 
$$
Furthermore, our assumptions are such that it is always possible 
to find a point $x^* \in S^1$ such that $u_{0x}(x^*)$ is negative. 
Setting $T_c=-1/u_{0x}(x^*)$ 
we conclude that the solution must blow-up in the $L^\infty$ norm 
as expected. 
\proofend 

More sophisticated break-down mechanisms for $\mu$DP can be demonstrated 
but we will not pursue them in this paper.


\subsection{Global solutions}

Our global result for $\mu$DP is obtained under a sign assumption 
on the initial data of the type that has been used previously 
in the studies of equations such as CH, $\mu$CH and DP. 
What follows is again a basic result. 
\begin{theorem}
Let $s>3$. Assume that $u_0 \in H^s(S^1)$ has non-zero mean and satisfies 
the condition
$$
\Lambda_\mu^2 u_0 \geq 0 
\quad
(\mathrm{or} \; \leq 0).
$$
Then the Cauchy problem for $\mu$DP has a unique global solution $u$ 
in $C(\mathbb{R},H^s(S^1)) \cap C^1(\mathbb{R},H^{s-1}(S^1))$.
\end{theorem}
\proofbegin From the local wellposedness and persistence result 
in Theorem \ref{th:lwp} 
we know that the solution $u$ is defined up to some time $T>0$ and that 
in order to extend it we only need to show that the norm 
$\| u_x(t, \cdot)\|_\infty$ remains bounded.

On the one hand, given any periodic function $w(x)$ differentiating 
the formula in \eqref{inverse} we readily obtain the estimate 
\begin{equation} \label{est} 
\| \partial_x w \|_{\infty} \lesssim \big\|\Lambda_\mu^2 w \big\|_{L^1}.
\end{equation} 

On the other hand, recall from Section \ref{Euler_Arnold} that 
for any solution $u$ of $\mu$DP with initial condition 
$u_0 \in H^s(S^1)$ the expression 
$m=\Lambda^2_\mu u$ is a tensor density in $\mathcal{F}_3(S^1)$ 
and hence from the $\mathrm{Diff}^s(S^1)$-action formula in \eqref{3-action} 
we get the following (pointwise) conservation law 
\begin{equation}  \label{local-mu} 
\frac{d}{dt} \Big( (\Lambda_\mu^2 u) \circ\xi \, (\partial_x \xi)^3 \Big)
=
\big(-u_{txx}-uu_{xxx}+3\mu(u)u_x -3u_xu_{xx}\big)\circ\xi \,(\partial_x \xi)^3
= 0, 
\end{equation} 
where $t \to \xi(t) \in \mathrm{Diff}^s(S^1)$ is the associated flow 
(i.e. $\dot{\xi} = u\circ\xi$ and $\xi(0)=e$). 
Thus, for any $x \in S^1$ as long as the solution exists 
it must satisfy 
$$
\Lambda_\mu^2 u\big(t,\xi(t,x)\big) \big(\partial_x \xi (t,x) \big)^3 
= 
\Lambda_\mu^2 u_0(x)
$$
and therefore by the assumption on $u_0$ it follows that
$$
\Lambda_\mu^2 u(t,x) \geq 0 
\quad 
\mathrm{for~any}~ x \in S^1 
$$
as long as it is defined.
The above inequality together with the conservation of the mean and \eqref{est} gives 
$$
\| \partial_x u(t, \cdot )\|_\infty
\lesssim
\big\|\Lambda_\mu^2 u(t) \big\|_{L^1}
= 
\int_0^1 \Lambda_\mu^2 u(t,x) \,  dx
= 
\mu(u_0) < \infty
$$
and we conclude that the solution $u$ must persist indefinitely.
\proofend


\subsection{Local wellposedness for any $\lambda$}
\label{subsec:lwp} 

For completeness we state here a local wellposedness result for 
the family (\ref{eq:mulambda}) of $\mu$-equations for any $\lambda$. 
Using the substitution $v=Av=\Lambda_\mu^2 u$ and the fact that, 
as in the case of $\mu$DP, the mean of the solution $u$ is conserved, 
we can write the Cauchy problem for \eqref{eq:mulambda} in the form 
\begin{equation} \label{L1}
u_t + u u_x 
+ 
\Lambda^{-2}_{\mu} \partial_x \left(
\lambda \mu(u) u + \frac{3-\lambda}{2} u_x^2  
\right)=0 
\end{equation}
\begin{equation} \label{L1-ic} 
u(0)=u_0. 
\end{equation} 
\begin{theorem}[Local wellposedness] \label{mueqlocalwellposedth}
Let $s>3/2$. For any $u_0\in H^s(S^1)$ there exist a $T > 0$ 
and a unique solution 
$$
u\in C\big( (-T,T),H^s \big) \cap C^1\big( (-T,T),H^{s-1}\big) 
$$
of \eqref{L1}-\eqref{L1-ic} which depends continuously 
on the initial data $u_0$. 
Furthermore, the solution persists as long as 
$\|u(t, \cdot )\|_{C^1}$ 
stays bounded. 
\end{theorem}
The proof of this theorem is analogous to Theorem \ref{th:lwp} and we omit it. 

In the previous section we proved a global existence result for $\mu$DP. 
Here, we formulate a similar result for any positive $\lambda$. 

For a fixed $\lambda >0$ let $u = u(t,x)$ be any smooth solution of 
the problem \eqref{L1}-\eqref{L1-ic} 
and let $\xi(t) \in \mathrm{Diff}^s(S^1)$ denote the corresponding flow 
of diffeomorphisms of the circle. 
As before, the argument rests on the pointwise conservation law 
obeyed by the solutions. Namely, using the formula for 
the action of $\mathrm{Diff}^s(S^1)$ on $\lambda$-densities 
$\mathcal{F}_\lambda(S^1)$ in \eqref{3-action} 
we find the analogue of \eqref{local-mu} to be 
$$
\frac{d}{dt}\Big( (\Lambda_\mu^2 u) \circ\xi \, (\partial_x \xi)^{\lambda} \Big)
=
\big( 
\Lambda^2_{\mu}u_t -uu_{xxx} +\lambda\mu(u)u_x -\lambda u_x u_{xx} 
\big)\circ\xi \, 
(\partial_x \xi)^{\lambda}
= 0
$$
so that as long as $u(t,x)$ exists 
$$
\Lambda_\mu^2 u\big( t,\xi(t,x)\big)\big(\partial_x\xi(t,x)\big)^{\lambda} 
= 
\Lambda_\mu^2 u_0(x), 
\qquad 
x \in S^1
$$
and therefore we have 
\begin{theorem}
Fix $\lambda>0$. Assume that $u_0 \in H^s(S^1)$ with $s>3$ has non-zero mean 
and satisfies the condition
$$
\Lambda_\mu^2 u_0 \geq 0
\quad
(\mathrm{or} \; \leq 0).
$$
Then the Cauchy problem \eqref{L1}-\eqref{L1-ic} has a unique global solution 
in $C\big(\mathbb{R},H^s(S^1)\big) \cap C^1\big(\mathbb{R},H^{s-1}(S^1)\big)$.
\end{theorem}
%


\section{Peakons}\nequation
\label{Waves} 
The CH and DP equations famously exhibit peakon solutions, see e.g. \cite{C-H, dhh}. 
In this section we show that each of the $\mu$-equations (\ref{eq:mulambda}) 
parametrized by $\lambda \in \Z$ also admits peakons. 
Let us point out here that the fact that $\mu$CH admits peaked solutions 
went unnoticed in \cite{K-L-M}, so that the present discussion provides 
an extension of the results of that paper. 

We will investigate existence of traveling-wave solutions of the form 
$u(t, x) = \varphi(x - ct)$. 
This will lead us to expressions for the one-peakon solutions. 
Subsequently, we will analyze the more general case of multi-peakons.


\subsection{Traveling waves}
We first consider equation (\ref{eq:mulambda}) where $\mu(u)$ is replaced by 
a parameter $\nu \in \R$, i.e.
\begin{equation}\label{lambdamuequation}
  -u_{txx} + \lambda \nu  u_x = \lambda u_x u_{xx} + uu_{xxx}.
\end{equation}
After analyzing traveling-wave solutions of this equation we will impose 
the two conditions 
$$
\mu(u) = \nu 
\qquad 
\mathrm{and}
\qquad  
\mathrm{period}(u) = 1
$$ 
to find which of these are traveling waves of equation (\ref{eq:mulambda}).
We assume that $\lambda \neq 0,1$. 
Substituting $u(t, x) = \varphi(x - ct)$ into equation (\ref{lambdamuequation}) 
we find 
$$
c\varphi_{xxx} + \lambda \nu \varphi_x 
= 
\lambda \varphi_x \varphi_{xx} + \varphi \varphi_{xxx}.
$$
We write this in terms of $f := \varphi - c$ and integrate to find
\begin{equation}\label{falphamu}  
  \lambda \nu f + a = \frac{\lambda - 1}{2} f_x^2 + ff_{xx},
\end{equation}
where $a$ in an integration constant.
We multiply equation (\ref{falphamu}) by $2 f^{\lambda - 2} f_x$ and integrate the resulting equation with respect to $x$. The outcome is
$$2\nu f^\lambda + \frac{2a}{\lambda - 1}f^{\lambda - 1} + b = f^{\lambda -1}f_x^2,$$
where $b$ is another constant of integration. This equation can be written as
\begin{equation}\label{varphiODE}  
  \varphi_x^2 = 2\nu(\varphi - c) + \frac{2a}{\lambda - 1} + \frac{b}{(\varphi - c)^{\lambda - 1}}.
\end{equation}
A solution of this equation yields the shape of a traveling wave within an interval where $\varphi$ is smooth, and gluing such smooth wave segments together at points where $\varphi = c$ produces the complete collection of traveling wave solutions of (\ref{lambdamuequation}). The analysis proceeds along the same lines and yields similar results as for CH (or as for HS if $\nu = 0$) cf. \cite{Ltrav}. Here we choose to consider the peaked solutions.

For completeness we point out that the equation analogous to (\ref{varphiODE}) when $\lambda = 0$ or $\lambda = 1$ is
\begin{equation} 
  \varphi_x^2 = \begin{cases} -2a + b(\varphi - c), & \lambda = 0, \\
  2\nu(\varphi - c) + 2a\mathrm{ln} |\varphi - c| + b, & \lambda = 1.
  \end{cases}
\end{equation}


\subsection{Periodic one-peakons}
The periodic peakons arise when $b = 0$ in (\ref{varphiODE}).\footnote{We henceforth assume $\nu \neq 0$.} Indeed, setting $b = 0$ and replacing $a$ by a new parameter $m$ via $a = (\lambda - 1) (c -m) \nu$ in (\ref{varphiODE}), we find
$$\varphi_x^2 = 2 \nu (\varphi - m).$$
Solving this equation, we find the peakon solution $u(t, x) = \varphi(x - ct)$ where
\begin{equation}\label{varphipeak}  
  \varphi(x) = m + \frac{\nu}{2} x^2 \quad \mathrm{for} \quad x \in \left[-\sqrt{\frac{2(c-m)}{\nu}}, \sqrt{\frac{2(c-m)}{\nu}}\right],
\end{equation}
and $\varphi$ is extended periodically to the real axis. The parameters are required to satisfy $(c-m)/\nu \geq 0$. The period of $\varphi$ is fixed by the condition that $\varphi = c$ at the peak. 

Equation (\ref{varphipeak}) defines a periodic peaked solution of (\ref{lambdamuequation}) for any choices of the real parameters $m, c, \nu$. 
In order to determine which of these solutions are solutions with period one of the $\mu$-equation (\ref{eq:mulambda}) we have to impose the conditions $\mathrm{mean}(\varphi) = \nu$ and $\mathrm{period}(\varphi) = 1$. A computation shows that $\mathrm{mean}(\varphi) = m + \frac{c - m}{3}$ and $\mathrm{period}(\varphi) = 2\sqrt{\frac{2(c-m)}{\nu}}$.
Therefore, $u(t, x) = \varphi(x -ct)$ where $\varphi$ is given by (\ref{varphipeak}) is a period-one solution of equation (\ref{eq:mulambda}) if and only if
$$2\sqrt{\frac{2(c-m)}{\nu}} = 1, \quad \frac{c-m}{\nu} \geq 0, \quad \mathrm{and} \quad \nu = m + \frac{c - m}{3}.$$
Solving these equations we find $m = 23c/26$ and $\nu = 12 c/ 13$. This leads to the following result.

\begin{theorem}
  For any $c \in \R$ and $\lambda \neq 0,1$, equation (\ref{eq:mulambda}) admits the peaked period-one traveling-wave solution $u(t, x) = \varphi(x - ct)$ where 
\begin{equation}\label{varphipeakon}  
  \varphi(x) = \frac{c}{26}(12 x^2 + 23)
\end{equation}  
  for $x \in [-\frac{1}{2}, \frac{1}{2}]$ and $\varphi$ is extended periodically to the real line. 
\end{theorem}
We note that the one-peakon solutions of (\ref{eq:mulambda}) are the same for any $\lambda$, that they travel with a speed equal to their height, and that $\mu(\varphi) = 12c/13$.


\subsection{Green's function}
For the construction of multi-peakons it is convenient to rewrite the inverse 
of the operator $\Lambda_\mu^2 = \mu - \partial_x^2$ in terms of a Green's function 
(an explicit formula for the inverse $\Lambda_\mu^{-2}$ was previously given in equation (\ref{inverse})).
We find
$$
(\Lambda_\mu^{-2}m)(x) = \int_0^1 g(x - x') m(x') dx',
$$
where the Green's function $g(x)$ is given by
\begin{equation}\label{greensdef}
g(x) 
= 
\frac{1}{2}x(x - 1) + \frac{13}{12} 
\quad 
\mathrm{for} 
\quad 
x \in [0, 1) \simeq S^1,
\end{equation}
and is extended periodically to the real line. 
In other words,
$$
g(x - x') 
= 
\frac{(x - x')^2}{2} - \frac{|x - x'|}{2} + \frac{13}{12}, 
\qquad 
x, x' \in [0, 1) \simeq S^1.
$$
Note that $g(x)$ has the shape of a one-peakon (\ref{varphipeakon}) 
with $c = 13/12$ and peak located at $x = 0$. In particular, $\mu(g) = 1$.


\subsection{Multi-peakons}
The construction of multi-peakons for the family of $\mu$-equations (\ref{eq:mulambda}) is similar to the corresponding construction for the family (\ref{CHDPfamily}) cf. \cite{dhh2}. Just like for CH and DP, the momentum $m = \Lambda_\mu^2u$ of an $N$-peakon solution is of the form
\begin{equation}\label{mpeakon}  
  m(t, x) = \sum_{i = 1}^N p_i(t) \delta(x - q^i(t)),
\end{equation}
where the variables $\{q^i, p_i\}_1^N$ evolve according to a finite-dimensional Hamiltonian system.
Indeed, assuming that $m$ is of the form (\ref{mpeakon}), the corresponding $u$ is given by
$$u = \Lambda_\mu^{-2}m = \sum_{i = 1}^N p_i(t) g(x - q^i(t)).$$ 
We assign the value zero to the otherwise undetermined derivative $g'(0)$, so that
\begin{equation}\label{greenderivative}
g'(x) := \begin{cases} 0, & x = 0, \\
x- \frac{1}{2}, & 0< x < 1.
\end{cases}
\end{equation}
This definition provides a naive way to give meaning to the term $-\lambda u_x m$ in the PDE (\ref{EF3}).
Substituting the multi-peakon Ansatz (\ref{mpeakon}) into the $\mu$-equation (\ref{eq:mulambda}) then shows that $p_i$ and $q^i$ evolve according to 
$$\dot{q}^i = \sum_{j = 1}^N p_j g(q^i - q^j), \qquad \dot{p}_i = - (\lambda - 1) \sum_{j = 1}^N p_i p_j g'(q^i - q^j).$$
A more careful analysis reveals that this computation is in fact legitimate and yields the following result. Note that the peakons fall outside the result of Theorem \ref{mueqlocalwellposedth} since for a peakon $u \notin H^s(S^1)$ for $s > 3/2$. 

\begin{theorem}\label{multipeakth}
The multi-peakon (\ref{mpeakon}) satisfies the $\mu$-equation (\ref{eq:mulambda}) in the weak form (\ref{L1})
in distributional sense if and only if $\{q^i, p_i\}_1^N$ evolve according to 
\begin{align}\label{peakonsystem}
  & \dot{q}^i = u(q^i), 
 	\\ \nonumber
  & \dot{p}_i = - (\lambda - 1)p_i \{u_x(q^i)\}
\end{align}  
where $\{u_x(q^i)\}$ denotes the regularized value of $u_x$ at $q^i$ defined by
\begin{equation}\label{uxpeakregularized}
  \{u_x(q^i)\} := \sum_{j = 1}^N  p_j g'(q^i - q^j)
\end{equation}
and $g'(x)$ is defined by (\ref{greenderivative}).
\end{theorem}
\proofbegin
  See Appendix \ref{multipeakapp}.
\proofend

\begin{figure}
\begin{center}
    \includegraphics[width=.3\textwidth]{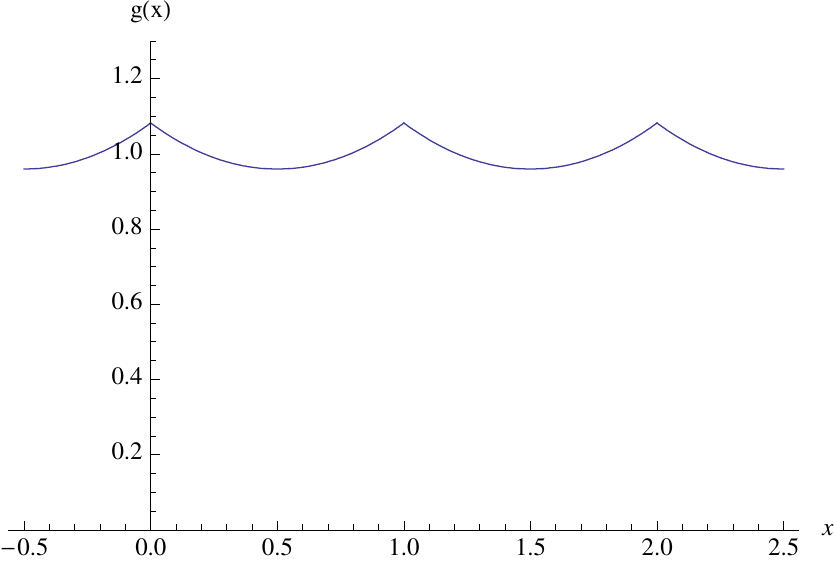} \quad
   \includegraphics[width=.3\textwidth]{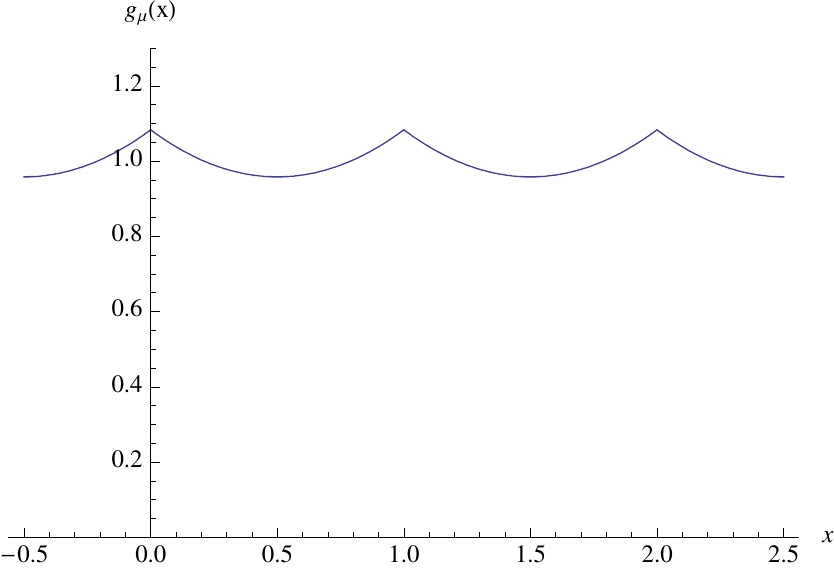}\quad
   \includegraphics[width=.3\textwidth]{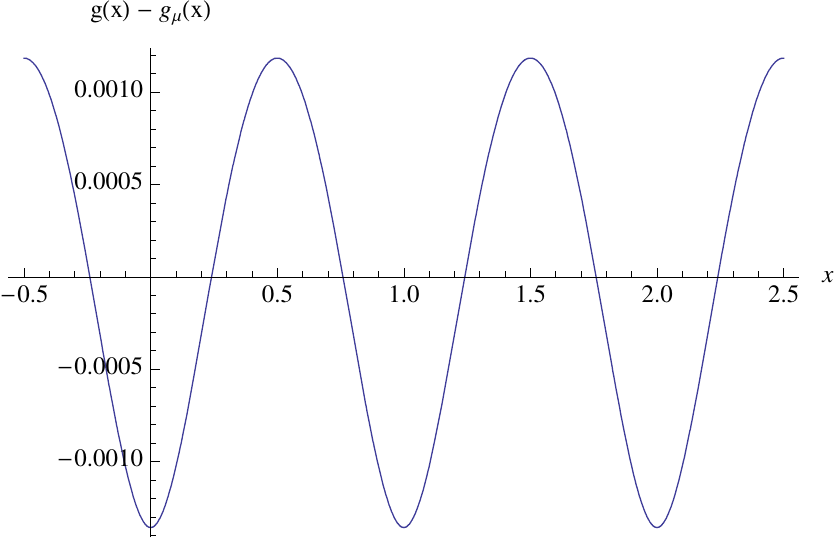} \\
     \begin{figuretext}\label{greensfigure} \upshape
       The periodic Green's functions $g(x)$ and $g_\mu(x)$ corresponding to the operators $A = 1 - \partial_x^2$ and $A = \mu - \partial_x^2$, respectively, and their difference $g(x) - g_\mu(x)$. 
     \end{figuretext}
     \end{center}
\end{figure}

The formulas for the peakons of the $\mu$-equations (\ref{eq:mulambda}) along with the expressions for the previously known peakons of the CH and DP family (\ref{CHDPfamily}) are summarized in Table \ref{peaktable}. At any particular time, the multi-peakon is a sum of Green's functions for the associated operator $A$, see Figure \ref{greensfigure}. Since the inertia operator for the $\mu$-equations ($A = \mu - \partial_x^2$) is different from that of the family (\ref{CHDPfamily}) ($A = 1 - \partial_x^2$), the new class of ``$\mu$-peakons'' have a different form with respect to the CH and DP peakons.
For a given family of equations the Green's function is the same for all $\lambda$, but the time evolution of the positions $q^i$ and momenta $p_i$ of the peaks of course depends on $\lambda$ and is for both families given by the system (\ref{peakonsystem}).

\begin{table}[htdp]
\begin{center}
\begin{tabular}{c|c|c}
Equation family	&	 Green's function ($x \in S^1 \simeq [0,1))$ & Multi-peakon $u(t,x)$
	\\ \hline
(\ref{CHDPfamily})	&
$g(x) = \frac{\hbox{cosh}(x - 1/2)}{2 \hbox{sinh(1/2)}}$
	&
$\sum_{i = 1}^N p_i(t) g(x - q^i(t))$
	\\
	& & 
	\\ 
(\ref{eq:mulambda})
	&
$g(x) = \frac{1}{2}\left(x - \frac{1}{2}\right)^2 + \frac{23}{24}$
	&
$\sum_{i = 1}^N p_i(t) g(x - q^i(t))$
	\\ 
\end{tabular}
\end{center}
\bigskip
\label{peaktable}
\caption{The periodic multi-peakon solutions of the family (\ref{CHDPfamily}) (which includes CH and DP) and of
the corresponding $\mu$-family (\ref{eq:mulambda}) (which includes $\mu$CH and $\mu$DP) when $\lambda \neq 0,1$.} 
\end{table}%

\subsection{Poisson structure}
The equations (\ref{peakonsystem}) take the Hamiltonian form
$$\dot{q}^i = \{q^i, H_0\}, \qquad \dot{p}_i = \{p_i, H_0\},$$
with respect to the Hamiltonian
$$H_0 = -\frac{\lambda^2}{\lambda -1} \sum_{i = 1}^N p_i,$$
and the Poisson structure
\begin{align*}
& \{p_i, p_j\} = \left(\frac{\lambda-1}{\lambda}\right)^2 p_i p_j G''(q^i - q^j), \qquad 
\{q^i, p_j\} = -\frac{\lambda-1}{\lambda^2} p_j G'(q^i - q^j),
	\\
& \{q^i, q^j\} = -\frac{1}{\lambda^2} G(q^i - q^j),
\end{align*}
where $G(x) = \int_0^x g(x')dx'$ cf. \cite{dhh2}. In order to make sense of $G(q^i - q^j)$ although $G$ is not periodic, we restrict $q^i, q^j$ to be local coordinates on $S^1$ satisfying say $0 < q^i, q^j < 1$. Then $-1 < q^i - q^j < 1$ and $G(q^i - q^j)$ is well-defined.
This Poisson structure can be derived from the Poisson structure corresponding to $J_0$ in (\ref{J0lambdadef}) by noting that
\begin{align*}
 \{m(x), m(y)\}_{J_0} =& -\frac{1}{\lambda^2} \bigl[G(x-y) m_x(x)m_x(y) - \lambda G'(x-y)m_x(x)m(y) 
 	\\
&+ \lambda G'(x-y)m(x) m_x(y) - \lambda^2 G''(x-y)m(x)m(y)\bigr].
\end{align*}
Indeed, let $m_0(x) = \sum_i p_{0i} \delta(x - q_0^i)$ be a multi-peakon of the form (\ref{mpeakon}). Let $\phi_i:S^1 \to \R$ (resp. $\psi_i:S^1 \to \R$) be a function which takes the value $1$ (resp. $x$)
in a small neighborhood of $x = q_0^i$ and the value $0$ elsewhere.
Then, for $m(x) = \sum_i p_{i} \delta(x - q^i)$ sufficiently close to $m_0$, we obtain
$$\int_{S^1} \phi_i(x) m(x) dx = p_i, \qquad \int_{S^1} \psi_i(x) m(x) dx = q^i p_i.$$
We now find $\{p_i, p_j\}$ from the computation
$$
\{p_i, p_j\} 
= 
\int_{S^1} \int_{S^1} \phi_i(x) \phi_j(y) \{m(x), m(y)\}_{J_0} dx dy 
= 
\left(\frac{\lambda-1}{\lambda}\right)^2 p_i p_j G''(q^i - q^j).
$$
We then compute
\begin{align*}
\{q^ip_i, p_j\} 
= &
\int_{S^1} \int_{S^1} \psi_i(x) \phi_j(y) \{m(x), m(y)\}_{J_0} dx dy 
	\\
=& -\frac{\lambda - 1}{\lambda^2} p_i p_j G'(q^i - q^j) + \left(\frac{\lambda-1}{\lambda}\right)^2 q^i p_i p_j G''(q^i - q^j).
\end{align*}
Since $\{q^ip_i, p_j\} = q^i\{p_i, p_j\} + \{q^i, p_j\}p_i$ and we know the value of $\{p_i, p_j\}$, this gives the expression for $\{q^i, p_j\}$. Finally, the expression for $\{q^i, q^j\}$ is obtained by calculating
$$\{q^ip_i, q^jp_j\} 
= 
\int_{S^1} \int_{S^1} \psi_i(x) \psi_j(y) \{m(x), m(y)\}_{J_0} dx dy$$
and comparing the result with the expression for
$$\{q^ip_i, q^jp_j\} = q^iq^j\{p_i, p_j\} + q^i\{p_i, q^j\} p_j + q^j\{q^i, p_j\}p_i + \{q^i, q^j\}p_ip_j$$
obtained by substituting in the expressions for $\{p_i, p_j\}$, $\{p_i, q^j\}$, and $\{q^i, p_j\}$.


\subsection{Multi-peakons for $\mu$CH}
For $\mu$CH ($\lambda = 2$) the equations (\ref{peakonsystem}) for $\{q^i, p_i\}$ are canonically Hamiltonian with respect to the Hamiltonian
$$h = \frac{1}{2}\sum_{i, j = 1}^N p_i p_j g(q^i - q^j),$$
and describe geodesic flow on $T^N \times \R^N$ (here $T^N = S^1 \times \cdots \times S^1$ denotes the $N$-torus) with respect to the metric $g_{ab}$ with inverse given by
$$g^{ij} = g(q^i - q^j).$$
In this case the reduction to multi-peakon solutions can be elegantly understood in terms of momentum maps using the ideas of \cite{H-M}. Indeed, let $S = \{1, 2, \dots, N\}$ be a set consisting of $N$ elements. We consider $(q^i, p_i)$ as an element of $T^*\hbox{Emb}(S, S^1) \simeq T^N \times \R^N$. Then the map
$$\mathbf{J}_{sing}: T^*\hbox{Emb}(S, S^1) \to \mathrm{vect}(S^1)^*:(q^i, p_i) \mapsto \sum_{i = 1}^N p_i \delta(x - q^i)$$
is an equivariant momentum map for the natural action of $\mathrm{Diff}(S^1)$ on $T^*\hbox{Emb}(S, S^1)$. Since momentum maps are Poisson, $\mathbf{J}_{sing}$ maps the solution curves in $T^*\hbox{Emb}(S, S^1)$ which satisfy (\ref{peakonsystem}) to solution curves in $\mathrm{vect}(S^1)^*$ which satisfy (\ref{CHtriple}). Note that $h = \mathbf{J}_{sing}^*H$ is the pull-back of the Hamiltonian $H$ on $\mathrm{vect}(S^1)^*$ defined by 
$H = \frac{1}{2} \int_0^1 um dx$.


\subsubsection{Periodic two-peakon for $\mu$CH}
In the case of two peaks for the $\mu$CH equation ($N = 2$, $\lambda = 2$), the system (\ref{peakonsystem}) is given explicitly by
\begin{align*}
 & \dot{q}^1 = p_1 g(0) + p_2 g(q^1 - q^2), \\
 & \dot{q}^2 = p_1 g(q^2 - q^1) + p_2 g(0), \\
 & \dot{p}_1 = - \dot{p}_2 = -p_1p_2 g'(q^1 - q^2).
 \end{align*}
Letting $Q = q^2 - q^1$ and $P = p_2 - p_1$ and using that 
$$h = \frac{13}{24}\frac{H_0^2}{16} + \frac{p_1p_2}{2}(Q^2 - |Q|),$$ 
we find the system
\begin{align*}
& \dot{Q} = -\frac{1}{2} P(Q^2 - |Q|),
	\\
& \dot{P} = -2\frac{\alpha}{Q^2 - |Q|}\left(Q - \frac{1}{2}\mathrm{sgn}(Q)\right),
\end{align*}
where $\alpha = 2h - \frac{13}{12}\frac{H_0^2}{16}$.
This leads to the following equation for $Q$:
$$
\ddot{Q}(Q^2 - |Q|) 
= 
(2\dot{Q}^2 + \alpha (Q^2 - |Q|))\left(Q - \frac{1}{2}\mathrm{sgn}(Q)\right).
$$
The solution to this ODE can be expressed in terms of inverses of elliptic functions.


\subsection{Shock-peakons for $\mu$DP}

In addition to the peakon solutions the DP equation admits an even weaker class 
of solutions with jump discontinuities called ``shock-peakons'', see \cite{lu}.
We will show here that $\mu$DP also allows shock-peakon solutions. 

We will seek the solutions of the form 
$$
m(t, x) 
= 
\sum_{i=1}^N \Big( 
p_i(t) \delta(x - q^i(t)) + s_i(t) \delta'(x - q^i(t))
\Big).
$$
The corresponding $u$ is
\begin{equation}\label{shockpeaku}
u = \sum_{i=1}^N \Big( p_i g(x - q^i) + s_i g'(x - q^i) \Big).
\end{equation}
Substituting this into the equation $m_t = -um_x - \lambda u_x m,$ we find formally
\begin{align*}
& \dot{q}^i = u(q^i), \\
& \dot{p}_i = (\lambda - 1) \left( s_i u_{xx}(q^i) - p_i u_x(q^i) \right), 
	\\
& \dot{s}_i = - (\lambda - 2)s_i u_x(q^i).
\end{align*}

The following theorem states that in the case of $\lambda = 3$ 
(this is the case of $\mu$DP), this formal computation can be 
rigorously justified provided that the $u_x$ and $u_{xx}$ terms 
are appropriately regularized. The reason that such weak solutions can be made sense of in the case of $\mu$DP is that the term $\frac{3-\lambda}{2} u_x^2$ in the weak formulation (\ref{L1}) is absent when $\lambda =3$.

\begin{theorem}\label{shockth}
The shock-peakon (\ref{shockpeaku}) satisfies the weak form (\ref{muDP}) of $\mu$DP 
in distributional sense if and only if $\{q^i, p_i, s_i\}_1^N$ evolve according to
\begin{align}\nonumber
& \dot{q}^i = u(q^i), 
	\\ \label{shocksystem}
& \dot{p}_i = 2 \big( s_i \{u_{xx}(q^i)\} - p_i \{u_x(q^i)\} \big),
	\\ \nonumber
& \dot{s}_i = - s_i \{u_x(q^i)\},
\end{align}
where
$$
\{u_x(q^i)\} 
= 
\sum_{j = 1}^N p_j g'(q^i - q^j) + \sum_{j = 1}^N s_j,
\qquad 
\{u_{xx}(q^i)\} = \sum_{j= 1}^N p_j,
$$
and $g'(x)$ is defined by (\ref{greenderivative}).
\end{theorem}
\proofbegin
See Appendix \ref{shockapp}.
\proofend


\section{The $\mu$Burgers equation} \nequation
\label{sec:3HS}

In this section we discuss the $\mu$B equation and its properties. 
We introduced this equation in Section \ref{Euler_Arnold} 
as the $\lambda=3$ analogue of the Hunter-Saxton equation. 
Here we take a different view and employ Riemannian geometric techniques 
to display its close relationship 
with the inviscid Burgers equation.


\subsection{Burgers' equation and the $L^2$-geometry of 
$\mathrm{Diff}^s \to \mathrm{Diff}^s/S^1$}
The constructions in this subsection are well known. 
First, recall that the inviscid Burgers equation 
\begin{equation} \label{B} \tag{B} 
u_t + uu_x =0 
\end{equation} 
is related to the geometry of the weak Riemannian metric on $\mathrm{Diff}^s(S^1)$ 
which on the tangent space $T_\xi\mathrm{Diff}^s(S^1)$ at a diffeomorphism $\xi$ 
is given by the $L^2$ inner product 
\begin{equation} \label{L2metric} 
\langle V, W \rangle_{L^2} 
= 
\int_0^1 V(x) W(x) \, dx 
\end{equation} 
where $V, W \in T_\xi\mathrm{Diff}^s$. 
It is not difficult to show that this metric admits a unique Levi-Civita connection 
$\nabla$ whose geodesics $\eta(t)$ in $\mathrm{Diff}^s(S^1)$ satisfy 
the equation\footnote{In fact, the geodesic equation is readily obtained from 
the first variation formula for the energy functional 
$E(\eta) = \frac{1}{2}\int_a^b \|\dot{\eta}(t) \|_{L^2}^2 dt$ 
of the $L^2$ metric \eqref{L2metric} on $\mathrm{Diff}^s(S^1)$.} 
\begin{equation} \label{B_geo} 
\nabla_{\dot\eta}\dot{\eta} = \ddot{\eta} 
= 
\big( u_t + u u_x \big)\circ\eta 
= 0 
\end{equation}
and hence correspond to (classical) solutions of the Burgers equation. 
Here $\eta(t)$ is simply the flow of $u(t,x)$ so that 
$\dot\eta(t,x) = u(t,\eta(t,x))$ 
and hence the second equality in \eqref{B_geo} follows at once from the chain rule. 

Consider next the homogeneous space $\mathrm{Diff}^s_0 = \mathrm{Diff}^s/S^1$ 
where $S^1$ denotes the subgroup of rotations. 
For $s>3/2$ this space is also a smooth Hilbert manifold whose tangent space 
at the point $[e]$ in $\mathrm{Diff}^s_0(S^1)$ corresponding to the identity 
diffeomorphism can be identified with the periodic Sobolev functions 
with zero mean 
\begin{equation*} \label{TDiff_0} 
T_{[e]}\mathrm{Diff}_0 = H^s_0(S^1) 
=
\left\{ w \in H^s(S^1): \int_0^1 w(x) \, dx =0  \right\}. 
\end{equation*}
It is useful to think of $\mathrm{Diff}^s_0(S^1)$ as the space of probability 
densities on the circle and view the map 
$\pi: \mathrm{Diff}^s \to \mathrm{Diff}^s_0$ 
as a submersion given by the pull-back 
$[\xi]=\pi(\xi) = \xi^\ast(dx)$ 
where $dx$ is the (fixed) probability measure with density 1 
(see e.g. \cite{ky}). 
The fibre through $\xi \in \mathrm{Diff}^s(S^1)$ consists of all
diffeomorphisms
that are obtained from $\xi$ by a composition on the left with a
rotation
and thus is just a right coset of the rotation subgroup $S^1$ which acts
on $\mathrm{Diff}^s(S^1)$ by left translations.
Furthermore, the $L^2$ metric \eqref{L2metric} is preserved by this action and 
$\pi$ becomes a Riemannian submersion with each tangent space 
decomposing into a horizontal and a vertical subspace 
$$
T_\xi\mathrm{Diff}^s 
= 
P_\xi(T_\xi\mathrm{Diff}^s) \oplus_{L^2} Q_\xi(T_\xi\mathrm{Diff}^s)
$$
which are orthogonal with respect to \eqref{L2metric}. 
The two orthogonal projections 
$P_\xi : T_\xi\mathrm{Diff}^s \to T_{\pi(\xi)}\mathrm{Diff}^s_0$ 
and 
$Q_\xi : T_\xi\mathrm{Diff}^s \to \mathbb{R}$ 
are given explicitly by the formulas 
\begin{equation}  \label{Pr} 
P_\xi(W) = W - \int_0^1 W(x) \, dx 
\qquad
\mathrm{and} 
\qquad
Q_\xi(W) = \int_0^1 W(x) \, dx. 
\end{equation} 

Finally, as with any Riemannian submersion\footnote{See e.g. O'Neill \cite{on} 
for details on Riemannian submersions.}, note that 
a necessary and sufficient condition for a curve $\eta(t)$ 
in $\mathrm{Diff}^s(S^1)$ to be an $L^2$ geodesic satifying \eqref{B_geo} 
(and hence correspond to a solution of Burgers' equation) 
is that 
$$
P_\eta \nabla_{\dot{\eta}}\dot{\eta} = Q_\eta\nabla_{\dot\eta}\dot{\eta}=0. 
$$

We are now ready to introduce the $\mu$B equation in this set-up.


\subsection{Lifespan of solutions, Hamiltonian structure and conserved quantities of the $\mu$Burgers equation} 

We first characterise local (in time) smooth solutions. 
\begin{theorem} \label{thm:3HS} 
A smooth function $u=u(t,x)$ is a solution of the $\mu$B equation 
$$
u_{txx} + 3u_x u_{xx} + u u_{xxx} = 0
$$
if and only if the horizontal component of the acceleration of 
the associated flow $\eta(t)$ in $\mathrm{Diff}^s(S^1)$ is zero 
i.e. 
$
P_\eta \nabla_{\dot\eta}\dot{\eta} = 0. 
$
In fact, given any $u_0 \in H^s(S^1)$ the flow of $u$ has the form 
$\eta(t,x) = x + t\big( u_0(x) - u_0(0) \big) + \eta(t,0)$ 
for all sufficiently small $t$. 
\end{theorem} 
\proofbegin 
Integrating the $\mu$B equation with respect to the $x$-variable we get 
\begin{equation} \label{3HS-int}
u_{tx} + u_x^2 + uu_{xx} = \int_0^1\big( u_x^2 + uu_{xx} \big)  dx = 0 
\end{equation} 
and integrating once again gives 
\begin{equation*}  
u_t + u u_x = c(t)
\end{equation*} 
where $c(t)$ is a function of the time variable only. 
Since the mean of a $\mu$B solution need not be preserved in time 
we obtain the equation\footnote{In the case when the mean $\mu(u)$ is 
independent of time \eqref{3HSB} becomes the standard (inviscid) Burgers equation. 
The presence of the mean is also the reason for our terminology.} 
\begin{equation} \label{3HSB} 
u_t + uu_x = \mu(u_t). 
\end{equation}
On the other hand, if $\eta(t)$ is the flow of $u$ in $\mathrm{Diff}^s(S^1)$, 
that is 
$$
\dot{\eta}(t,x) = u(t,\eta(t,x)), 
\qquad 
\eta(0,x) = x 
$$
then from \eqref{B_geo} and \eqref{Pr} we compute the horizontal component 
of the acceleration of $\eta(t)$ in $\mathrm{Diff}^s(S^1)$ to be 
\begin{align*} 
P_\eta \nabla_{\dot\eta}\dot{\eta} 
= 
\ddot{\eta} - \int_0^1 \ddot\eta \, dx 
= 
\big( u_t + uu_x \big)\circ\eta 
- 
\int_0^1 \big(u_t + uu_x \big)\circ\eta \, dx. 
\end{align*} 
Using these formulas and integrating by parts we find that 
the equation 
$P_\eta\nabla_{\dot\eta}\dot\eta = 0$ 
is equivalent to 
\begin{align}  \label{calc} 
\big( u_t + uu_x \big)\circ\eta 
= 
\ddot\eta 
= 
\int_0^1 \ddot\eta \, dx 
&= 
\int_0^1 \ddot\eta\circ\eta^{-1} \, dx    \\ \nonumber
&= 
\int_0^1 ( u_t + u u_x ) \, dx 
= 
\int_0^1 u_t \, dx 
\end{align}
which establishes the first part of the theorem. 

In order to prove the second statement it suffices to observe that 
from \eqref{calc} we have in particular 
$$
\ddot\eta (t,x) 
= 
\int_0^1 \ddot{\eta}\circ\eta^{-1}(t,x) \, dx 
=
\ddot\eta(t,0) 
$$
which immediately implies 
$$
\eta(t,x) - \eta(t,0) = x+ t\big( u_0(x) - u_0(0) \big) 
$$
for any $0 \leq x \leq 1$ and any $t$ for which $\eta$ is defined. 
\proofend


\begin{remark} 
It is easy to verify that equation \eqref{3HSB} is also equivalent to 
$\nabla_{\dot\eta} P_\eta \dot{\eta} = 0$ 
which can be interpreted as saying that the horizontal component of the velocity 
of the flow $\eta(t)$ of the $\mu$B equation is parallel transported along 
the flow. 
\end{remark} 

\begin{remark} 
A similar construction can be used to derive the $\mu$B equation 
from the \emph{right-invariant} $L^2$ metric on $\mathrm{Diff}^s(S^1)$ 
defined by 
$$
\langle V, W \rangle_{\xi} 
= 
\int_0^1 V\circ\xi^{-1}(x) W\circ\xi^{-1}(x) \, dx 
\quad 
\mathrm{where} 
\;\;
V, W \in T_\xi \mathrm{Diff}^s 
\;\; 
\mathrm{and} 
\;\;
\xi \in \mathrm{Diff}^s. 
$$ 
In this case the corresponding orthogonal projections are 
$
\tilde{P}_\xi(W) = W - \int_0^1 W \, d\xi = W - \tilde Q_\xi(W) 
$
and the geodesic equation in $\mathrm{Diff}^s(S^1)$ reads
$
\tilde{\nabla}_{\dot\eta}\dot{\eta}=\ddot\eta+2 \dot{\eta} \, 
\partial_x\dot{\eta} \, (\partial_x\eta)^{-1}=0.
$
Proceeding as in the proof of Theorem \ref{thm:3HS} we then get 
$$
0 = \tilde{P}_\eta \tilde{\nabla}_{\dot\eta}\dot\eta  
= 
\big( u_t + 3u u_x - \mu(u_t) \big)\circ\eta 
$$ 
which after rescaling the dependent variable yields \eqref{3HSB}. 
\end{remark} 


It is therefore not surprising to find that the $\mu$B equation shares 
a number of properties with the Burgers equation. 
\begin{corollary} \label{thm:3HSB} 
Suppose that $u(t,x)$ is a smooth solution of the $\mu$B equation 
and let $u(0,x)=u_0(x)$. 
\begin{enumerate} 
\item 
The following integrals are conserved by the flow of $u$ 
$$ 
\int_0^1 \big( u-\mu(u) \big)^p dx 
= 
\int_0^1 \big( u_0-\mu(u_0) \big)^p dx, 
\qquad
p=1,2,3 \dots.
$$ 
\item 
There exists $T_c >0$ such that 
$\| u_x(t)\|_\infty \nearrow \infty$ as $t \nearrow T_c$.
\item
The $\mu$B equation has a bi-Hamiltonian structure given by the Hamiltonian functionals
$$
H_0=-\frac{9}{2}\int m \ dx, \ \ H_2=-\frac{1}{6}\int u^3 dx, 
$$
and the operators
$$
J_0= -m^{2/3}\partial_x m^{1/3}\partial_x^{-3}m^{1/3}\partial_x m^{2/3}, \ \ J_2= \partial_x^5.
$$
\end{enumerate} 
\end{corollary} 
\proofbegin 
The first statement follows by direct calculation. 
Regarding the second, it suffices to consider the equation in \eqref{3HS-int} 
and apply the argument that was used in the proof of Theorem \ref{muDPblowup}. 
The third statement follows from a straightforward computation establishing
$$
J_0 \frac{\delta H_0}{\delta m}=J_2 \frac{\delta H_2}{\delta m}=-um_x-3u_x m
$$
where $m=-\partial_x^2 u$ and our earlier observation from section \ref{Hamilton} that $J_0$ and $J_2$ are compatible.
\proofend


\section{Multidimensional $\mu$CH and $\mu$DP equations}\nequation
\label{sec:EPDiff}

In this final section we briefly consider possible generalizations of 
the $\mu$CH and $\mu$DP equations to higher dimensions.
Let $\mathfrak{X}(T^n)$ be the space of smooth vector fields on 
the $n$-dimensional torus $T^n$ comprising the Lie algebra of 
$\mathrm{Diff}(T^n)$. 
Let $\mathbf{A}$ be a self-adjoint positive-definite operator defining 
an inner product on the Lie algebra. 
Given a vector field $\mathbf{u} \in \mathfrak{X}(T^n)$ we define 
the corresponding momentum density by $\mathbf{m} = \mathbf{Au}$.
The EPDiff equation for geodesic flow on $\mathrm{Diff}(T^n)$ is given by 
\begin{equation}\label{EPdiff}
\mathbf{m}_t + \mathcal{L}_\mathbf{u} \mathbf{m} = 0, 
\end{equation}
where $\mathcal{L}_\mathbf{u} \mathbf{m}$ denotes the Lie derivative of the momentum one-form density $\mathbf{m}$ in the direction of $\mathbf{u}$ (see \cite{H-M}). 
Identifying $T^n \simeq \R^n/\Z^n$ and letting $x^i$, $i = 1, 2, \dots, n$ 
denote standard coordinates on $\R^n$, this can be written as
$$
\frac{\partial m_i}{\partial t} 
+ 
u^j \frac{\partial m_i}{\partial x^j} 
+ 
m_j \frac{\partial u^j}{\partial x^i} + m_i \text{div}(\mathbf{u}) 
= 0 
$$
where $\mathbf{m} = m_i dx^i \otimes d^nx$ 
and 
$\mathbf{u} = u^i \partial/\partial x^i$.
The CH equation is the $(1+1)$-dimensional version of 
the EPDiff equation (\ref{EPdiff}) when $\mathbf{A} = 1 -\Delta$. 
It is natural to define 
a multidimensional generalization of $\mu$CH as equation (\ref{EPdiff}) 
with the operator $\mathbf{A}$ defined by
$$
\mathbf{A}\mathbf{u} 
= 
(\mu - \Delta)\mathbf{u} 
= 
\int_{T^n} \mathbf{u} \, d^nx - \Delta \mathbf{u}.
$$


Similarly, using higher-order tensor densities we arrive at 
multidimensional versions of DP and $\mu$DP. 
Let $\mathbf{m} = m_i dx^i \otimes d^n x \otimes d^n x$. 
Then equation (\ref{EPdiff}) becomes
\begin{equation}\label{EPdiffDP}
\frac{\partial m_i}{\partial t} 
+ 
u^j \frac{\partial m_i}{\partial x^j}  
+ 
m_j \frac{\partial u^j}{\partial x^i} + 2 m_i \text{div}(\mathbf{u}) 
= 0.
\end{equation}
In $1+1$ dimensions this reduces to
$$
m_t + m_x u + 3m u_x = 0.
$$
Hence, when $\mathbf{A} = 1 - \Delta$, equation (\ref{EPdiffDP}) 
can be viewed as a multidimensional DP equation, 
while if $\mathbf{A} = \mu - \Delta$ it can be viewed as 
a multidimensional $\mu$DP equation. 


\appendix
\section{Proof of Theorem \ref{multipeakth} }\label{multipeakapp}
\renewcommand{\theequation}{A.\arabic{equation}}\nequation
In this appendix we prove Theorem \ref{multipeakth}, that is, we show that the multi-peakons defined in (\ref{mpeakon}) are weak solutions in the distributional sense of the $\mu$-equation (\ref{eq:mulambda}) if and only if $\{q^i, p_i\}_1^N$ evolve according to (\ref{peakonsystem}).

Equation (\ref{eq:mulambda}) in weak form reads\footnote{Here $*$ denotes convolution: for two functions $f,g:S^1 \simeq [0,1) \to \R$, we have $(f * g)(x) = \int_0^1 f(x - y)g(y) dy$.} 
\begin{equation}\label{mueqgreenform}  
  u_t + \frac{1}{2}(u^2)_x + \lambda \mu(u) g' * u + \frac{3-\lambda}{2} g' * (u_x^2)= 0,
\end{equation}
where $g$ is the Green's function defined in (\ref{greensdef}).
This is to be satisfied in the space of distributions $\mathcal{D}'(\R \times S^1)$, i.e.
$$\int_{-\infty}^\infty \int_{S^1} \left( u\phi_t + \frac{1}{2} u^2 \phi_x + \lambda \mu(u)(g * u)\phi_x + \frac{3 - \lambda}{2}(g * u_x^2)\phi_x \right)dxdt = 0,$$ 
for all test functions $\phi(t, x) \in C_c^\infty(\R \times S^1)$. 

Let $g_j := g(x - q^j)$ and $g_j' := g'(x - q^j)$. We specify the value of $g'$ at $x = 0$ by setting $g'(0) := 0$ as in equation (\ref{greenderivative}).
We will need the following identities which can be verified 
by direct computation
\begin{align} \label{gconv1}
  g' * g_j = &-\frac{1}{3}\left(g_j - \frac{13}{12} \right)g_j',
		\\ \label{gconv2}
  g' * g_j' =& 1 - g_j,
		\\ \label{gidentity1}
-(g_i - g_j)(g_i' - g_j') =& 2g_j'(q^i)g_i + 2g_i'(q^j) g_j
	\\ \nonumber
& + \left(g(q^i - q^j) - \frac{13}{12}\right)g_i' + \left(g(q^j - q^i) - \frac{13}{12}\right)g_j', 
		\\ \label{gidentity2}
g_i' g_j' =& g_i + g_j + g'(q^i - q^j)(g_i' - g_j') + g(q^i - q^j) - 3.
\end{align}
These identities hold pointwise for all $x, q^i, q^j \in S^1$ (except (\ref{gidentity2}) which holds for all $x, q^i, q^j \in S^1$ unless $x = q^i = q^j$; a fact that does not matter in what follows).
Let $u = \sum_i p_i g_i$ be a peakon. As explained in the appendix of \cite{lu}, 
we can compute the left-hand side of (\ref{mueqgreenform}) as a distribution 
in the variable $x$ without having to involve test functions explicitly. 
This yields
$$
u_t = \sum_{i=1}^N \left(\dot{p}_i g_i - p_i \dot{q}^i g_i'\right),
\qquad 
(u^2)_x = \sum_{i,j =1}^N 2 p_ip_j g_ig_j'.
$$
Using (\ref{gconv1}) we find
$$
\lambda \mu(u) g' * u 
= 
-\frac{\lambda}{3} 
\sum_{i,j=1}^N p_i p_j\left(g_j - \frac{13}{12}\right)g_j' .
$$
Moreover,
$$
u_x^2 = \sum_{i,j =1}^N p_ip_j g_i'g_j' ,
$$
which in view of (\ref{gconv1}), (\ref{gconv2}), and (\ref{gidentity2}) 
implies that
$$
g' * (u_x^2) 
= 
-\frac{1}{3} \sum_{i,j =1}^N p_ip_j \left(
(g_ig_i' + g_jg_j' - \frac{13}{12}(g_i' + g_j') 
+ 
3 g'(q^i - q^j)(g_i - g_j) \right).
$$
Using these ingredients we can write (\ref{mueqgreenform}) as
$$
\sum_{i=1}^N \left( \dot{p}_i g_i - p_i \dot{q}^i g_i' \right) 
+ 
\sum_{i,j =1}^N p_ip_j g_ig_j'
- 
\frac{\lambda}{3} \sum_{i, j=1}^N p_i p_j\left(g_j - \frac{13}{12}\right)g_j'
$$
$$
- \frac{3-\lambda}{6}  \sum_{i,j =1}^N p_ip_j 
\left( (g_ig_i' + g_jg_j' - \frac{13}{12}(g_i' + g_j') 
+ 3 g'(q^i - q^j)(g_i - g_j) \right) 
= 0.
$$
We rewrite this as
$$
\sum_{i=1}^N \left\{ 
\dot{p}_i g_i + p_i\left(\frac{13}{12}
\left(\sum_j p_j\right) - \dot{q}^i\right) g_i' \right\} 
- 
\frac{1}{2} \sum_{i \neq j} p_ip_j (g_i - g_j)(g_i' - g_j')
$$
$$
+ \frac{\lambda - 3}{2}  \sum_{i,j =1}^N p_ip_j g'(q^i - q^j)(g_i - g_j) = 0.
$$
Using (\ref{gidentity1}) we find
\begin{align*}\label{sumdotpi2pi}
\sum_{i=1}^N \left\{ 
\left(\dot{p}_i +   2 p_i \sum_{j \neq i} p_j g_j'(q^i)\right) g_i
+ p_i\left(\frac{13}{12}p_i 
+ 
\sum_{j \neq i} p_jg(q^i - q^j) - \dot{q}^i\right) g_i' \right\} 
	\\ \nonumber
+ 
\frac{\lambda - 3}{2}  \sum_{i,j =1}^N p_ip_j g'(q^i - q^j)(g_i - g_j) = 0.
\end{align*}
The last sum on the left-hand side can be written as
$$
(\lambda - 3)  \sum_{i,j =1}^N p_ip_j g' \big( q^i - q^j \big) g_i, 
$$
so that, using the definition (\ref{uxpeakregularized}) of $\{u_x(q^i)\}$ 
and the fact that $g(0) = 13/12$, we arrive at
$$
\sum_{i=1}^N \left( (\dot{p}_i +  (\lambda - 1) p_i\{u_x(q^i)\}) g_i
+ p_i(u(q^i) - \dot{q}^i) g_i' \right) = 0.
$$
Since $g_i$ and $g_i'$ form a linearly independent set, this equation holds 
if and only if $\{q^i, p_i\}_1^N$ evolve according to (\ref{peakonsystem}). 
This completes the proof of Theorem \ref{multipeakth}.


\section{Proof of Theorem \ref{shockth}}\label{shockapp}
\renewcommand{\theequation}{B.\arabic{equation}}\nequation
In this appendix we prove Theorem \ref{shockth}, that is, we show that 
the shock-peakons defined in (\ref{shockpeaku}) are weak solutions 
in the distributional sense of $\mu$DP if and only if $\{q^i, p_i, s_i\}_1^N$ 
evolve according to (\ref{shocksystem}). The objective of the proof is similar 
to the corresponding proof for DP presented in the appendix of \cite{lu}. 
However, since we consider the spatially periodic case and the Green's functions 
are very different, the details of the two proofs are quite different.

The $\mu$DP equation in weak form reads
\begin{equation}\label{muDPgreens}
u_t + \frac{1}{2}(u^2)_x + 3\mu(u) g' * u = 0,
\end{equation}
which is to be satisfied in distributional sense.
We use the same notation as in the proof of Theorem \ref{multipeakth}. 
Let $\delta_i(x) = \delta(x - q^i)$ and note that $g_i'' = 1 - \delta_i$.
It is enough to verify that (\ref{muDPgreens}) holds as 
a distribution in $x$ when the functions $u_t$, $\frac{1}{2}(u^2)_x$, 
and $3\mu(u) g' * u$ are replaced by
\begin{align*}
& u_t = \sum_{i=1}^N \left( 
\dot{p}_i g_i - p_i \dot{q}^i g_i' + \dot{s}_i g_i' 
- s_i \dot{q}^i(1 - \delta_i)
\right),
	\\
& \frac{1}{2}(u^2)_x 
= 
\sum_{i,j =1}^N \left( p_ip_j g_ig_j' + s_is_j g_i'(1 - \delta_j) 
+ 
p_is_j (g_i' g_j' + g_i(1 - \delta_j))\right),
\end{align*}
and
\begin{equation}\label{3mugprimestar}
3\mu(u) g' * u 
= 
\sum_{i, j=1}^N \left( 
-p_i p_j\left(g_j - \frac{13}{12}\right)g_j'  + 3 p_i s_j (1 - g_j) 
\right),
\end{equation}
respectively.  
In order to derive (\ref{3mugprimestar}) we used the fact that
$$
g' * u = \sum_{j=1}^N \left( p_j(g' * g_j) + s_j (g' * g_j')\right)
$$
together with the identities (\ref{gconv1}) and (\ref{gconv2}).

Putting these expressions together we find that equation (\ref{muDPgreens}) 
can be written as
\begin{align*}
& \sum_{i=1}^N \left( 
\dot{p}_i g_i - p_i \dot{q}^i g_i' + \dot{s}_i g_i' - s_i \dot{q}^i(1 - \delta_i)
\right)
  	\\
 & + \sum_{i,j =1}^N 
\left( p_ip_j g_ig_j' +  s_is_j g_i'(1 - \delta_j) 
+ 
p_is_j (g_i' g_j' + g_i(1 - \delta_j))\right) 
 	\\
&+ \sum_{i, j=1}^N \left( 
-p_i p_j\left(g_j - \frac{13}{12}\right)g_j'  + 3 p_i s_j (1 - g_j)
\right) = 0.
\end{align*}
Employing (\ref{gidentity1}) we can rewrite this as
\begin{align} \nonumber
& \sum_{i=1}^N \left\{ 
\left(\dot{p}_i + 2 p_i \sum_{j \neq i} p_j g_j'(q^i) \right) g_i
+ 
p_i\left(\frac{13}{12}p_i 
+ 
\sum_{j \neq i} p_jg(q^i - q^j) - \dot{q}^i\right) g_i' 
+ \dot{s}_i g_i' - s_i \dot{q}^i\right\} 
	\\ \label{alongshockeq}
& + \sum_{i,j =1}^N \left( 
s_is_j g_i' + p_is_j (g_i' g_j' + g_i) 
\right)
	\\ \nonumber
& +3 \sum_{i, j=1}^N p_i s_j (1 - g_j) 
+ 
\sum_{i=1}^N s_i\left(\dot{q}^i - \sum_j s_j g_j' - \sum_j p_j g_j\right) \delta_i 
= 0.
\end{align}
The identity (\ref{gidentity2}) shows that the two terms 
$\sum_{i,j} p_is_j (g_i' g_j' + g_i)$ 
and 
$3 \sum_{i, j} p_i s_j (1 - g_j)$ 
combine to give
$$
\sum_{i,j =1}^N p_i s_j \left( 
2(g_i - g_j) + g'(q^i - q^j)(g_i' - g_j') + g(q^i - q^j)
\right).
$$
We find that (\ref{alongshockeq}) can be written as
\begin{align*}
& \sum_{i=1}^N \left(\dot{p}_i 
+ 
2 p_i \sum_{j \neq i} p_j g_j'(q^i)  
+ 
2p_i \sum_j s_j - 2s_i \sum_j p_j \right) g_i 
	\\
& + \sum_{i=1}^N \left(p_i\left(\sum_j p_jg_j(q^i) 
+ 
\sum_j s_jg'(q^i - q^j) - \dot{q}^i\right)  + \dot{s}_i 
+ 
s_i \sum_j s_j - s_i \sum_j p_j g'(q^j - q^i)\right) g_i' 
	\\
& + \sum_{i=1}^Ns_i\left(\sum_j p_jg(q^i - q^j) -  \dot{q}^i\right)
+ 
\sum_{i=1}^N s_i\left(\dot{q}^i - \sum_j s_j g_j' - \sum_j p_j g_j\right) \delta_i 
= 0.
\end{align*}
The term $\sum_i s_i\left(\sum_j p_jg(q^i - q^j) -  \dot{q}^i\right)$ equals
$$
\sum_i s_i\left(\sum_j p_jg(q^i - q^j) + \sum_j s_j g_j'(q^i) -  \dot{q}^i\right).
$$
Using that
$$
u(q^i) = \sum_j p_j g_j(q^i) + \sum_j s_j g_j'(q^i),
$$
$$
\{u_x(q^i)\} = \sum_j p_j g_j'(q^i) + \sum_j s_j, \qquad \{u_{xx}(q^i)\} = \sum_j p_j,
$$
we arrive at the equation
$$
\sum_{i=1}^N \left(\dot{p}_i + 2 p_i\{u_x(q^i)\} 
- 
2s_i \{u_{xx}(q^i)\} \right) g_i 
+ 
\sum_{i=1}^N \left(p_i\left(u(q^i)  - \dot{q}^i\right)  
+ 
\dot{s}_i + s_i \{u_x(q^i)\}\right) g_i' 
$$
$$ 
+ \sum_{i=1}^Ns_i(u(q^i) -  \dot{q}^i)
+ \sum_{i=1}^N s_i\left(\dot{q}^i - u(q^i)\right) \delta_i = 0.
$$
Since $\{1, g_i, g_i', \delta_i\}$ form a linearly independent set 
this equation holds if and only if $\{q^i, p_i, s_i\}_1^N$ 
evolve according to (\ref{shocksystem}). 
This completes the proof of Theorem \ref{shockth}.


\section{Conservation of $H_2$ for $\mu$DP}\label{H2app}
\renewcommand{\theequation}{C.\arabic{equation}}\nequation
In this appendix we verify explicitly that the functional $H_2$ defined in (\ref{H0H2muDPdef}) is conserved under the flow of $\mu$DP. The following two forms of $\mu$DP are used:
\begin{equation*}
u_t + u u_x +3\mu(u)\Lambda_\mu^{-2} u_x =0
\end{equation*}
and its derivative
\begin{equation*}
u_{tx}+uu_{xx}+u_x^2+3(\mu(u))^2-3u\mu(u)=0.
\end{equation*}
We also use the identities $\mu (\Lambda_\mu^{-2} u)=\mu(u)$ and $\Lambda_\mu^{-2} u_{xx}= -u+\mu(u)$.
We compute
\begin{align*}
\frac{d H_2}{dt} = &-\int_{S_1} \biggl(3\mu(u)(\Lambda_\mu^{-2} u_x)( \Lambda_\mu^{-2} u_{xt} )+\frac{1}{2}u^2 u_t\biggr)
	\\
= &\int_{S_1} \biggl( 3 \mu(u)( \Lambda_\mu^{-2}  u_x) \Lambda_\mu^{-2}\left(uu_{xx}+u_x^2+3\mu(u)^2-3u\mu(u)\right) 
	\\
& + \frac{1}{2}u^2\left(uu_x+3\mu(u)\Lambda_\mu^{-2} u_x \right) \biggr).
\end{align*}
We use $uu_{xx}+u_x^2=\frac{1}{2}\partial_x^2(u^2)$ and find
\begin{align*}
\frac{d H_2}{dt}  =& \int_{S_1} \biggl( -9\mu(u)^2 (\Lambda_\mu^{-2} u_x) \Lambda_\mu^{-2} u + 9 \mu(u)^3\Lambda_\mu^{-2} u_x + \frac{3}{2} \mu(u)(\Lambda_\mu^{-2}\partial_x^2(u^2))\Lambda_\mu^{-2} u_x
	\\
&+ \frac{1}{2}u^3u_x + \frac{3}{2} u^2 \mu(u)\Lambda_\mu^{-2} u_x \biggl).
\end{align*}
Here the first two terms and the fourth term vanish by periodicity. The third term is equal to
\begin{equation*}
  \frac{3}{2} \mu(u)(\Lambda_\mu^{-2}\partial_x^2(u^2))\Lambda_\mu^{-2}u_x= \frac{3}{2}\mu(u)(-u^2+\mu(u^2))\Lambda_\mu^{-2} u_x
\end{equation*}
whose first part cancels the fifth term and the second part vanishes by periodicity.


\bigskip
\noindent
{\bf Acknowledgement} {\it J.L. is grateful to Professor D. D. Holm for valuable discussions and suggestions. J.L. acknowledges support from a Marie Curie 
Intra-European Fellowship.}


\bibliography{is}

\end{document}